\documentclass[journal, one column, draftcls, 12pt]{IEEEtran}
\usepackage{hyperref}
\usepackage{graphics, subfigure,  times, amsfonts, amsmath, amssymb, cite}
\usepackage{tikz}
\usepackage{url}
\usetikzlibrary{shapes,arrows}
\usepackage{pgfplots}
\usepackage{color}
\usepackage{multirow}
\usetikzlibrary{automata,positioning}
\usepackage{rotating}
\usepackage{flushend}
\usepackage{amsthm}
\usepackage{epstopdf}
\newtheorem{theorem}{Theorem}
\newtheorem{assumption}{Assumption}
\newtheorem{lemma}{Lemma}
\newtheorem{definition}{Definition}

\usepackage{bbm}

\title{ Menu-Based Pricing for Charging of Electric Vehicles with Vehicle-to-Grid Service }
\date{}
\author{Arnob Ghosh and Vaneet Aggarwal\thanks{The authors are with the School of  Industrial Engineering, Purdue University, West Lafayette, IN 47907, email: \{ghosh39,vaneet\}@purdue.edu. This work was supported in part by the U.S. National Science Foundation
		under grant CCF-1527486.}}

\begin{document}
	\maketitle
\begin{abstract}
The paper considers a bidirectional power flow model of the electric vehicles (EVs) in a charging station. The EVs can inject energies by discharging via a Vehicle-to-Grid (V2G) service which can enhance the profits of the charging station. However, frequent charging and discharging degrade  battery life. A proper compensation needs to be paid to the users to participate in the V2G service. We propose a menu-based pricing scheme, where the charging station selects a price for each arriving user for the amount of battery utilization, the total energy, and the time (deadline) that the EV will stay. The user can accept one of the contracts or rejects all depending on their utilities. The charging station can serve users using a combination of the renewable energy and the conventional energy bought from the grid. We show that though there exists a profit maximizing price which maximizes the social welfare, it provides no surplus to the users if the charging station is aware of the utilities of the users. If the charging station is not aware of the exact utilities, the social welfare maximizing price may not maximize the expected profit. In fact, it can give a zero profit. We propose a pricing strategy which provides a guaranteed fixed profit to the charging station and it also maximizes the expected profit for a wide range of utility functions. Our analysis shows that when the harvested renewable energy is small the users have higher incentives for the V2G service. We, numerically, show that the charging station's profit and the user's surplus both increase as V2G service is efficiently utilized by the pricing mechanism. 
\end{abstract}
\vspace{-0.1in}
 \section{Introduction}
\subsection{Motivation}
Electric Vehicles (EVs) have several advantages over the traditional gasoline powered vehicles. For example, EVs are more environment friendly and more energy efficient. Realizing the above, regulators (e.g. Federal Energy Regulator Commission (FERC)) are providing incentives to the consumers to switch to electric vehicles. Manufacturers (e.g. Tesla, Nissan) are increasingly developing EVs equipped with superior technologies. As a result, electric vehicles are increasingly become popular.   However, a wide deployment of EVs requires an extensive network of charging stations which can be  capable of charging large number of vehicles. 

Vehicle-to-grid (V2G) service has been proposed \cite{Tomic,Kempton} to enhance the profitability of the EVs. In the V2G service, EVs can inject energies to the grid by discharging from their batteries. Thus, this bidirectional power flow where EVs can both charge and discharge has a lot of potential.  Hence, a lot of effort is going on for developing bidirectional EVs\cite{tesla}.  However, the higher amount of charging and discharging cycles will degrade the battery life. Hence, the owners of the EVs have to be compensated adequately for the V2G service. Thus, though a charging station can gain an additional profit using the V2G service of the EVs, however without a proper pricing mechanism, the owners of the EVs will not prefer the V2G service in the first place which may nullify the profit of the charging station. 

Without a proper control mechanism, the cost of the charging station and the peak energy consumption may increase. Without a profitable charging station, the wide-scale deployment of the EVs will remain a distant dream. The charging station needs to select prices in order to earn profits by allocating resources in an intelligent manner among the EVs. The charging station also needs to provide adequate compensation to the owners if the EV is used for the V2G service. However, high prices or low compensation may not provide incentives to the owners which may reduce the profit.  Hence, a proper pricing mechanism for charging the EVs and the V2G service is imperative for a charging station. 

\vspace{-.1in}
\subsection{Our Contributions}
We propose a menu based pricing scheme for charging an EV. Whenever an EV arrives at the charging station, the charging station offers  a variety of  contracts $(l,t,BU)$ to the EV's owner (or, user) at a price $p^{BU}_{l,t}$  to the user where the user will be able to charge at least $l$ units of energy within the deadline $t$  for completion, and the battery usage will be limited to $l+BU$ amount. The battery usage is the total amount of charging and discharging of the EV. The user either accepts one of the contracts by paying the specified price or rejects all of those based on its payoff. We assume that the user gets a utility for consuming $l$ amount of energy within the deadline $t$. However, the user also has to incur a cost for battery utilization $BU$. The payoff of the user (or, user's surplus) for a contract is the difference between the utility and the sum of the cost incurred, and the price to be paid for the contract. The user will select the option which fetches the highest payoff. The various advantages of the above pricing scheme should be noted. First, it is an online pricing scheme. It can be adapted for each arriving user.  Second, since the charging station offers prices for different levels of charging required and the deadline, the charging station can prioritize one contract over the others depending on the energy resources available.  Third, the charging station also provides options of the maximum battery utilization to the users. The user who is not interested in the V2G service is entitled to do that by selecting the contract with $BU$ as $0$.  Finally,  the user's decision is much simplified.  She only needs to select one of the contracts (or, reject all). 

We consider that the charging station is equipped with renewable energy harvesting devices and a storage device for storing energies. The charging station may also  buy conventional energy from the market to fulfill the contract of the user if required. Hence, if a new user accepts  the contract $(l,t,BU)$, a cost is incurred to the charging station. This cost may also depend on the existing EVs and their resource requirements. A contract also specifies the maximum battery utilization which restricts the V2G service generated from an EV. Hence, the charging station needs to find the optimal cost for each contract. We show that obtaining the cost of fulfilling a contract is a {\em linear programming} problem. We also show that if a user accepts a contract with a higher battery utilization, the cost of the contract is lower (Lemma 1). 
% it has the above favorable property,

We consider two optimization problems--i) social welfare\footnote{Social welfare is the sum of the profit of the charging station and the user surplus.} maximization, and  ii) the EV charging station's profit maximization. We investigate the existence of a pricing mechanism which maximizes the ex-post social welfare, {\em i.e.}, maximizes the social welfare  for every possible realization of the utility function.   We show that there exists such a pricing strategy. The pricing scheme is  simple to compute, as the charging station selects a price which is equal to the marginal cost for fulfilling a certain contract for a new user (Theorem~\ref{thm:max_socv2g}).   However, the above pricing scheme only provides {\em zero} profit to the charging stations. Thus, such a pricing scheme may not be useful to the charging station. We show that when a charging station is {\em clairvoyant} (i.e., the primary knows the utilities of the users), there exists a pricing scheme which satisfies both the objectives (Lemma~\ref{thm:profit_max}). Though  in the above pricing mechanism, the user's surplus becomes $0$. Thus, a {\em clairvoyant} charging station may not be beneficial for the user's surplus.

  In the scenario where the charging station does not know the exact utilities of the users, we show that {\em there may not exist a pricing strategy which simultaneously maximizes the ex-post social welfare and  the expected profit}. One has to give away the ex-post social welfare maximization in order to achieve expected profit maximization.   {\em However, the user's surplus becomes higher compared to the clairvoyant scenario}. Hence, an uncertainty of the utility enhances the user's surplus. We propose a pricing strategy which can fetch the highest possible profit to the charging station under the condition that it  maximizes the ex-post social welfare (Theorem~\ref{thm:profitmax_uncertainty}).   Above pricing strategy provides a {\em worst case} maximum profit to the charging station. % We show that the charging station's profit increases as the renewable energy generation increases.  % The user's surplus also shows a similar trend. 

 Since the above pricing strategy  may not yield the {\em maximum expected} profit to the charging station, we have to relax the constraint the social welfare to be maximized in order to yield a higher profit to the charging station. Whether a contract will be selected by the user does not depend on the price of the contract, but also, the prices of other contracts. Thus, achieving a pricing scheme which maximizes the expected profit is difficult because of the discontinuous nature of the profits. We propose a pricing strategy which yields a fixed (say $\beta$) amount of profit to the charging station.  Further, we show that a suitable choice of $\beta$ can maximize the profit of the charging station for a class of utility functions (Theorem~\ref{thm:aclassutility}). 

In Section~\ref{sec:v2g_profitability} we characterize the conditions which will yield higher profits to the charging station in the V2G service. We show that if the conventional energy price is high, the charging station selects more preferable prices for enticing V2G services if the renewable energy harvested is low. When the charging station's storage capacity is low and the renewable energy harvesting is low, the user's incentives for the V2G service also increases. The V2G service also increases the profit of the charging station. %{\em When the harvested renewable energy is very high, the users have lower incentives.} However, when the harvested energy is limited, the charging station selects more favorable prices to the users in order to incentivize for the V2G service. %However, if users participate in the V2G service, the users may have  However,  if users participate in the v2G service, the cost to serve a contract again decreases, and the users may have higher incentives to participate in the V2G service. 

 Finally, we, empirically provide insights how a trade-off between the profit of the charging station and the social welfare can be achieved for various pricing schemes (Section~\ref{sec:simulation_results}).  Our numerical analysis shows that both the user's surplus and the charging station's profit increase with the V2G service. The energy consumption during the peak period decreases as users provide more V2G services during the peak period. %We also show that because of the V2G service, lower amount of energies are consumed during the peak hours. The charging station also attains higher profits. %We also cha the amount of energy feed back from the vehicles. 

\subsection{Related Literature} Charging schedule for EVs using price signals for {\em unidirectional} service have been considered\cite{anciliary_service2,tong,sortomme,javidi}.  The above papers did not consider the optimal discharging schedule of the EVs as these papers only considered unidirectional power flow. As a result, these paper did not consider the battery degradation cost incurred in the V2G service. 

Few papers considered the {\em bidirectional} power flow  \cite{sortommev2g,wu_tsg,poor_tsg,wong,he,xing}. However, their focus was scheduling of the charging and discharging pattern of the EVs,  rather than the pricing aspects.  Naturally, these papers did not consider whether users will prefer the amount of battery degradation found in the optimal scheduling process. The social welfare maximizing and profit maximizing prices are also not considered.  \cite{schaar} considered the optimal pricing to the EVs in a day-ahead setting for residential charging. The users  control the charging and discharging pattern at each instance for the price selected by the aggregator in a pre-specified manner. However, the users can arrive randomly in the charging station, the charging station needs to select a price for each arriving user using an online algorithm. The users can not control the charging and discharging schedule at each instance in the charging station.  Compared to {\em all} the above papers, in our approach, the charging station specifies the amount of battery utilization for each contract. The charging station also selects prices for different deadlines in our approach. Hence, the user can now choose its own deadline and can specify the V2G service it is willing to commit based on its own need. %Compared to {\em all} the above papers, in our menu-based pricing approach, we consider controlling of the deadline, and the battery utilization of the EVs by selecting intelligent pricing mechanism. 

Deadline differentiated prices have been considered \cite{bitar,bitar2,nayyar}. In our previous work \cite{tech_ev}, we also considered a menu-based pricing scheme where the charging station provides prices for different deadlines and prices to the new user. However, the above papers did not consider the {\em bidirectional} power flow problem. Thus, in this setting the charging station now also needs to select different prices for different amounts of battery utilization. {\em To the best of our knowledge, this is the first attempt to consider menu-based pricing which considered the bidirectional power flow from the EVs.}  Further, this is the first work that  incentivizes the users to participate in the V2G service by finding optimal prices for different amounts of battery degradation while maximizing the social welfare or profit of the charging station.%The charging station needs to find the optimal prices for different amounts of battery degradation is not considered in the above papers. The user's additional cost is not considered. Naturally, the profit of the charging station for V2G service, the incentives for the users to participate in the V2G have been not considered in those papers. %However, those paper considered that they have some information about the vehicle's driving patterns which may be difficult to obtain, or the residential charging. However, in our approach, we consider an online pricing mechanism. 
 \begin{figure*}
 \begin{minipage}{0.49\linewidth}
 \includegraphics[trim=0in 0.2in 0in 0in,width=0.99\textwidth, height=40mm]{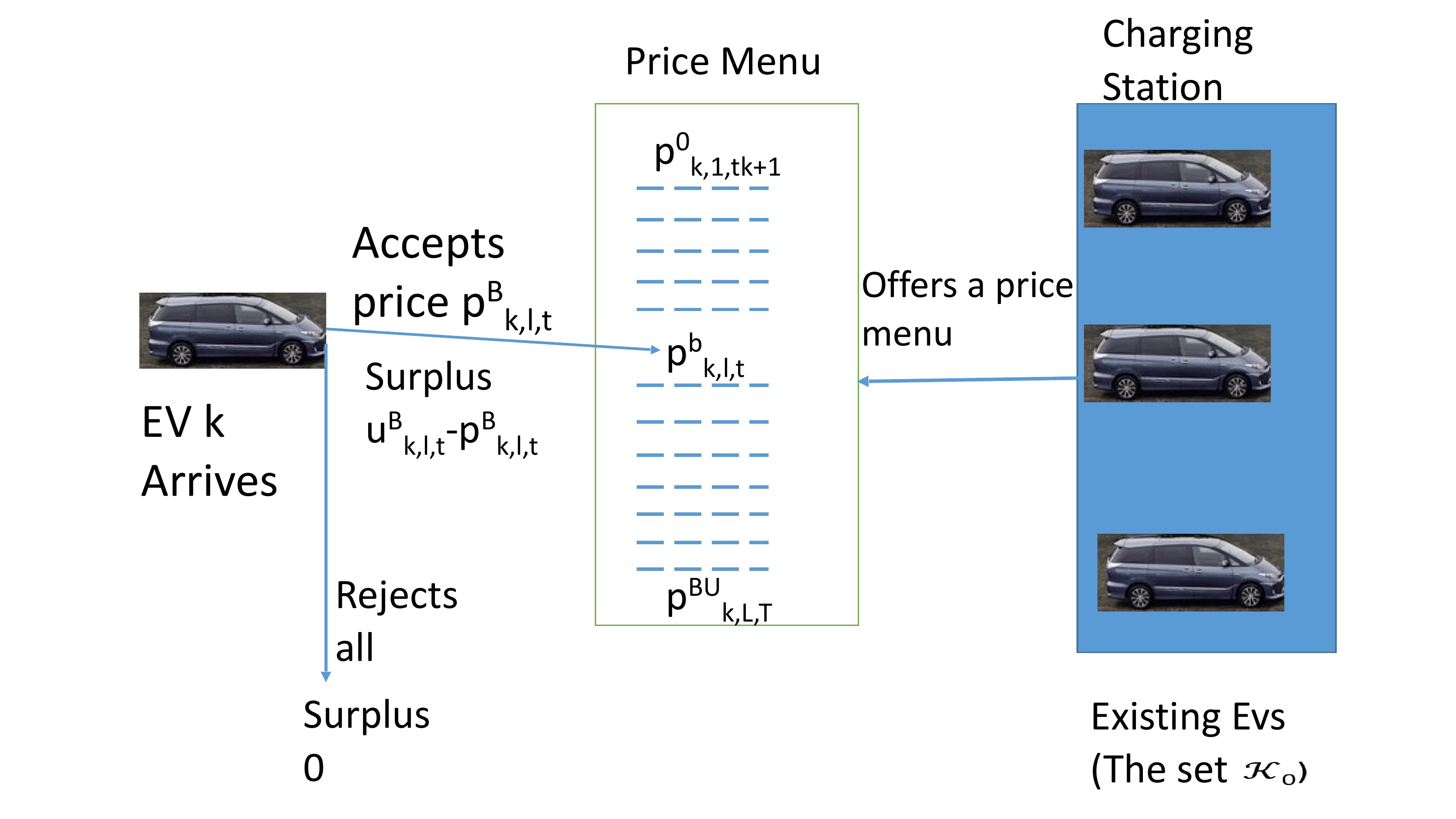}
 \vspace{-0.2in}
 \caption{The trading model: Charging station offers a menu of contracts for $l$, $t$, and $BU$; the arriving user decides either one of them or rejects all.}
 \label{fig:ev_charging}
 \vspace{-0.2in}
 \end{minipage}\hfill
 \begin{minipage}{0.49\linewidth}
 \includegraphics[trim=0in 0.7in 0in 0in,width=0.99\textwidth,height=38mm]{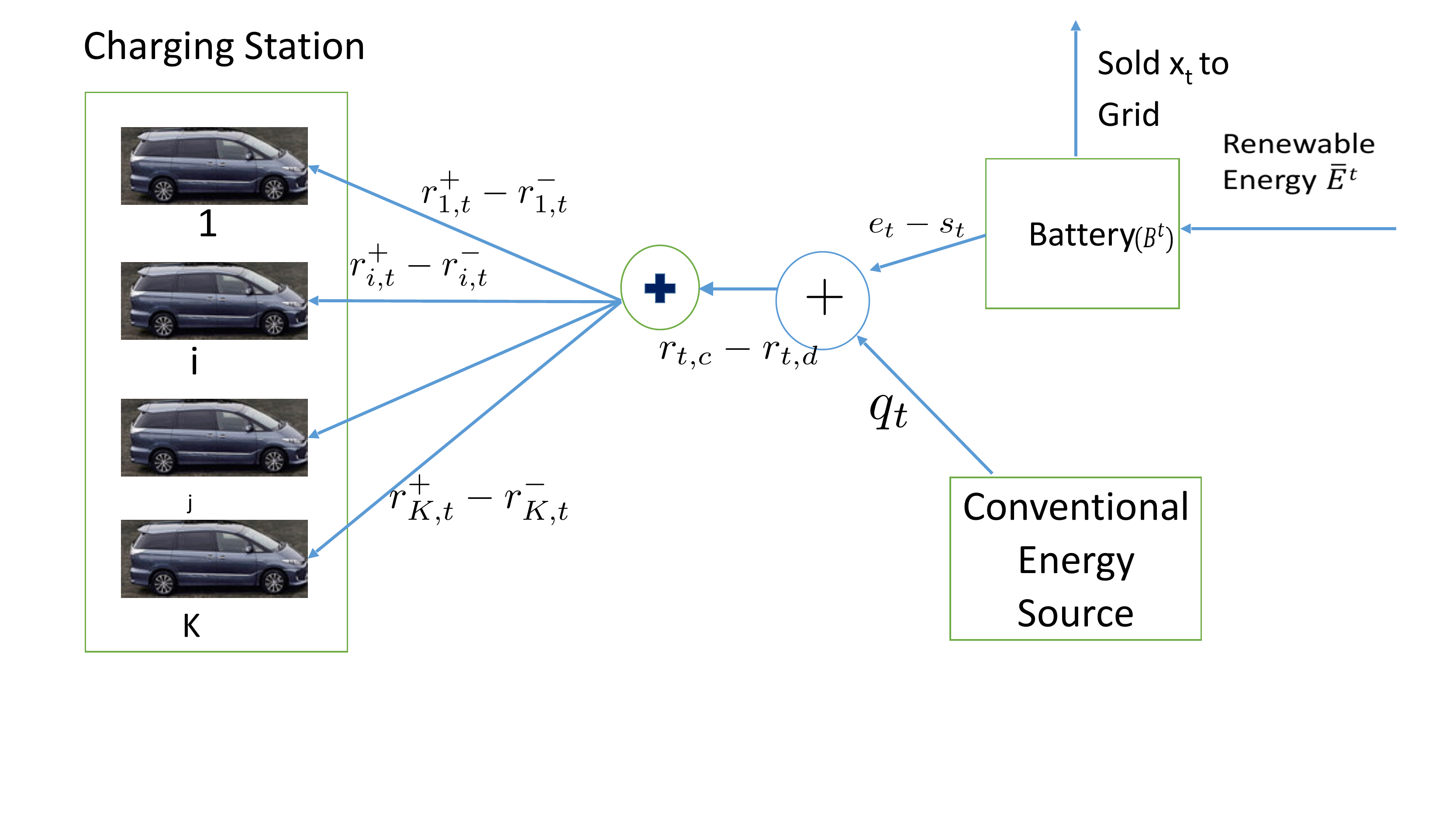}
 \vspace{-0.2in}
 \caption{The hybrid energy source, the battery, and the charging and discharging of EVs. $r_{t,c}, r_{t,d}$ denote $r_{t,\text{charge}}, r_{t,\text{discharge}}$ respectively.}
 \label{fig:hybrid_v2g}
 \vspace{-0.2in}
 \end{minipage}
 \end{figure*}
 \section{System Model}
 \subsection{Menu-based Pricing for arriving user}
 We consider that EVs arrive throughout a day at the charging station for charging.  Suppose that the user $k$ arrives at time $t_k$. The job of the charging station is to select a price for charging in order to maximize the profit. We consider a vehicle-to-grid (V2G) service where electric vehicles (EVs) can feed back  stored power  to the grid. Specifically, the energy can be discharged from the batteries of the EVs. However, the EV batteries have fixed number of charging-discharging cycles \cite{Delucchi,zhou}. When the battery is used for feeding back energy to the grid, the battery wear cost may be significant. Thus, the users may want to limit the battery utilization as low as possible. 
 
 We consider that a charging station wants to maximize its profit over a certain time period (e.g. over a day). It will offer a menu-based price contract $p^{BU}_{k,l,t}$  for contract $(l,t,BU)$ to the user $k$ which arrives at time $t_k$ in the charging station (Fig.~\ref{fig:ev_charging}). If the user selects the menu $(l,t,BU)$, then, the user $k$ will be able to charge $l$ amount of energy at most within the deadline $t$ and the maximum amount of {\em additional} battery utilization (Definition 1) is restricted to $BU$. The deadline $t$ denotes that the user has to leave by time $t$,  after which the EV will not be able to charge.
 %However, since the EVs may be bi-directional, thus, the charging station can sell the stored energies from the EVs and can earn some profits. However, the EVs also have to be compensated for discharging as frequent discharging and charging can lead to degradation of battery life. 
 
 \subsection{System Constraints}
First, we describe the constraints the charging station has to satisfy the contract $(l,t_{\text{dead}},BU)$ for the user $k$. 

\textbf{Charging/Discharging rate}: Let $r_{k,t}$ be the amount of energy provided to (or, discharged from) the EV $k$ during time $[t,t+1)$. $r_{k,t}>0$ indicates that the EV $k$ is charged and $r_{k,t}<0$ indicates that energy is discharged from the EV $k$ during time $[t,t+1)$.  Also note that there is an initial set $\mathcal{K}_0$ of existing EVs. Vehicle $i\in \mathcal{K}_0$ requires additional $N_i$ amount of energy within the deadline $w_i$. The charging and discharging efficiency of the EV $k$ is denoted as $\eta_{k,c}\leq 1$ and $\eta_{k,dc}\leq 1$ respectively.  Note that $r_{k,t}=r^{+}_{k,t}-r^{-}_{k,t}$ where $r^{+}_{k,t}$ is the positive part of $r_{k,t}$ (it denotes that the electric vehicle $k$ is charged during time $[t,t+1)$) and $r^{-}_{k,t}$ denotes the negative part of $r_{k,t}$ ( it denotes the amount discharged from EV $k$ during time $[t,t+1)$).  Hence, the following set of constraints must be satisfied--
\begin{eqnarray}
& \sum_{t=t_k}^{t_{\text{dead}}-1}r_{k,t}\geq l,\quad \sum_{t=t_k}^{w_i-t_k-1}r_{i,t}\geq N_{i}, \forall i\in \mathcal{K}_0\label{eq:charged}\\
& r_{t,\text{charge}}=r^{+}_{k,t}/\eta_{k,c}+\sum_{i\in \mathcal{K}_0}r^{+}_{i,t}/\eta_{i,c},\label{eq:charge}\\
& r_{t,\text{discharge}}=\eta_{k,dc}r_{k,t}^{-}+\sum_{i\in \mathcal{K}_0}\eta_{i,dc}r_{i,t}^{-}\label{eq:chargeanddischarge}\\
& r_{k,t}=r_{k,t}^{+}-r_{k,t}^{-},\quad r_{i,t}=r_{i,t}^{+}-r_{i,t}^{-}\quad \forall i\in \mathcal{K}_0.\label{eq:chargeeachev}
\end{eqnarray}
 $r_{t,\text{charge}}$ indicates the total amount of energy required for charging  and $r_{t,\text{discharge}}$ indicates the total amount of energy used when EVs discharge during time $[t,t+1)$.  Note that an EV can not simultaneously charge and discharge, hence, we must have $r^{+}_{i,t}r^{-}_{i,t}=0$ for all $i\in \mathcal{K}_0\cup \{k\}$ and $t$. {\em Later we will show in Theorem~\ref{thm:noposinegative} that in an optimal solution $r^{+}_{i,t}r^{-}_{i,t}=0$ though we have not explicitly considered the above constraint in the system.}

\textbf{EV's battery limit}: Let the battery level of the EV $i$ at time $t$ be $\text{EV}_{i,t}$. Let the battery level at the start of the time $t_k$ be $\text{EV}^{i,{\rm ini}}$. Since the battery of the EV can be charged or discharged in the V2G service, the total amount of charging and discharging must satisfy the limits on the battery levels:The battery level must be between a lower value $d_{i,\text{min}}$ (it can be $0$\footnote{However, low depth of discharge can degrade the battery life\cite{zhou}. Hence, the minimum value can be set at some positive value.} ) and the high value $d_{i,{\rm max}}$ (the highest capacity). Hence,  we must have
\begin{align}\label{eq:dischargelimit}
&\text{EV}_{i,t+1}=\text{EV}_{i,t}+r^{+}_{i,t}-r^{-}_{i,t},\nonumber\\
& \text{EV}_{i,t_k}=\text{EV}^{i,{\rm ini}}\quad  \forall i\in \mathcal{K}^{0}\cup \{k\},\nonumber\\
& d_{i,{\rm min}}\leq \text{EV}_{i}\leq d_{i,{\rm max}} \quad \forall i\in \mathcal{K}^{0}\cup \{k\}.
\end{align}
Note that $r^{+}_{i,t}-r^{-i}_{i,t}$ is $r_{i,t}$.

\textbf{Charging and Discharging rate limit}:
There is also a charging and discharging limit. Hence, we have 
\begin{eqnarray}\label{eq:limit}
R_{min}\leq r_{k,t}\leq R_{max},\quad R_{min}\leq r_{i,t}\leq R_{max} \forall i\in \mathcal{K}_0.
\end{eqnarray}
$R_{max}$, and $R_{min}$ are respectively the charging and discharging rate limits. 

\textbf{Hybrid Source}: The charging station is equipped with renewable energy harvesting devices (Fig.~\ref{fig:hybrid_v2g}). The charging station also has a storage device with capacity $B_{max}$. The charging station has the forecast of harvested energy as $\bar{E}^t$ for time $[t,t+1)$. The charging station also can buy energies from the conventional market.  Suppose that the amount of energy used from the storage device during time $[t,t+1)$ be $e_t$ and let the energy be stored from the electric vehicles and the conventional energy in the storage device be $s_t$ during time $[t,t+1)$ (Fig.~\ref{fig:hybrid_v2g}). Let $B^{t}$ be the level of battery at time $t$. The charging efficiency and discharge efficiency of the battery of the charging station is considered to be $\eta_{c,cs}$ and $\eta_{d,cs}$ respectively. Thus, 
 \begin{eqnarray}\label{eq:battery_capacity}
 & B^{t+1}=B^{t}+\bar{E}^t\eta_{c,cs}-e_t/\eta_{d,cs}+s_t\eta_{c,cs}, \label{eq:battery_level}\\
 & B_{max}\geq B^{t+1}\geq 0, B^{t_k-1}=B^{0}, B^{T}=B^{0}.\label{eq:battery_capacity}
 \end{eqnarray}
Note that our model can also incorporate the scenario where the batteries have some static leakage rate i.e. the battery level decreases with time.

 \textbf{Energy to and from the grid}: Let the amount the charging station buys from the conventional market  and sells to the grid as $q_t$ and $x_t$ respectively for time $[t,t+1)$ (Fig:~\ref{fig:hybrid_v2g}). Hence, we have
 \begin{align}
 e_t=\max\{r_{t,\text{charge}}-q_t+x_t-r_{t,\text{discharge}},0\}\nonumber\\
 \quad s_t=\max\{q_t-r_{t,\text{charge}}+r_{t,\text{discharge}}-x_t,0\}.\nonumber
 \end{align}
 Since $s_t$ is the negative of $e_t$, thus,  we can represent the above constraint in the following--%we can eliminate the variable $e_t$ and we denote $s_t\geq 0$ if the energy is stored and $s_t<0$ if the energy is discharged respectively from the storage device of the charging station during time $[t,t+1)$. We denote the above equality as
 \begin{align}\label{eq:storage}
 q_t-r_{t,\text{charge}}+r_{t,\text{discharge}}-x_t=s_t-e_t, s_t\geq 0, e_t\geq 0. \end{align}
 We also must have $e_ts_t=0$ since the battery of the charging station should not charge or discharge at the same time. Though we have not explicitly considered the above constraint, however, we show that in an optimal solution in Theorem~\ref{thm:noposinegative}, we always have $e_ts_t=0$. 
 
 Note that the constraints in (\ref{eq:battery_capacity}) and (\ref{eq:storage}) also specify the bound on the amount of charging and discharging from the EVs. The amount of charging can not exceed the total amount of stored energy in the battery and the energy bought from the grid. The stored energy depends on the harvested renewable energy, energy charged to the EVs and discharged from the EVs, energy bought and energy sold to the grid.
 
 Fig.~\ref{fig:hybrid_v2g} depicts various system parameters.%We assume that the estimated amount $\bar{E}^t$ matches the renewable energy generation. %amount that can be used for charging and discharging from the EVs is upper bounded by the total renewable energy generation, the 
 
 \textbf{Maximum Battery Utilization}
 In the vehicle-to-grid (V2G) service,  frequent charging and discharging may reduce the lifecycle of the battery of an EV. Thus, the users may not like the EVs be charged and discharged very often. We define a metric which will model the total maximum utilization of the battery of an EV. 
 \begin{definition}
 Battery utilization is defined as the absolute value of the difference between the battery levels  at two subsequent time intervals. 
 
 Hence, if the deadline is $t_{dead}$ for user $k$, then the total battery utilization for user $k$ is $\sum_{t=t_k}^{t_{\text{dead}}-1}|\text{EV}_{k,t+1}-\text{EV}_{k,t}|$.
 \end{definition}
%Note that similar metric is used for battery utilization in various papers \cite{he}.  
The battery utilization defines the total level of charging and discharging has been done. If the EV $k$ is used only for charging, then $\sum_{t=t_{k}}^{t_{\text{dead}}-1}|\text{EV}_{k,t+1}-\text{EV}_{k,t}|=l$. In a contract $(l,t_{dead},BU)$, the battery of EV needs a charging amount of $l$, and thus, the total battery utilization has to be at least $l$. %However, an user may not prefer higher battery utilization as it would increase the battery wear cost. 
 
 Thus, the contract $(l,t_{\text{dead}},BU)$ where $BU=0, 1,\ldots, BU_{max}$ denotes the {\em additional battery utilization} apart from the charging amount $l$ specified by the contract. We denote $BU$ as the maximum {\em additional} battery utilization with slight abuse of notation. Specifically,  in the contract $(l,t_{\text{dead}},BU)$, the maximum utilization  is restricted to $l+BU$  for user $k$. 
 
 Note that $l$ is the energy that the user $k$ will receive in the contract $(l,t_{\text{dead}},BU)$. Suppose the maximum utilization remaining for an existing user $i\in \mathcal{K}_0$ at time $t_k$ is $BU_i$. Thus, if the user $k$ selects the contract $(l,t_{\text{dead}},BU)$, the constraint that the charging station has to satisfy is
 \begin{align}\label{eq:bu}
 \sum_{t=t_k}^{t_{\text{dead}}-1}|\text{EV}_{k,t+1}-\text{EV}_{k,t}|-l\leq BU,\nonumber\\
 \quad \sum_{t=t_k}^{w_i-t_k-1}|r_{i,t}-r_{i,t-1}|-N_i\leq BU_i \forall i\in \mathcal{K}_0.
 \end{align}
 The above constraint is not linear. In the following we reduce it in a linear form. 
 
 Note from (\ref{eq:dischargelimit}) that 
 \begin{align}\label{eq:identity}
& \text{EV}_{i,t+1}=\text{EV}_{i,t}+r^{+}_{i,t}-r^{-}_{i,t}=r_{i,t}+\text{EV}_{i,t}.\nonumber\\
& |\text{EV}_{i,t+1}-\text{EV}^{i,t}|=|r_{i,t}|.
 \end{align}Since $r^{+}_{i,t}$ and $r^{-}_{i,t}$ are the positive and negative parts of $r_{i,t}$, thus, $|r_{i,t}|=r^{+}_{i,t}+r^{-}_{i,t}$.  Hence, the expression in (\ref{eq:bu}) becomes
 \begin{align}\label{eq:battery_utilization}
  \sum_{t=t_k}^{t_{\text{dead}}-1}r^{+}_{k,t}+r^{-}_{k,t}-l\leq BU,\nonumber\\ \sum_{t=t_k}^{w_i-t_k-1}r^{+}_{i,t}+r^{-}_{i,t}-N_i\leq BU_i \forall i\in \mathcal{K}_0.
\end{align}
Note that we must have $r^{+}_{i,t}r^{-}_{i,t}=0$. We show in Theorem~\ref{thm:noposinegative} that this is indeed true in an optimal solution.%\footnote{In \cite{he}, non-linear constraints have also been considered for battery utilization. Our model can easily be extended to those scenarios, however, then the problem will be non-linear convex optimization problem. The computational complexity will increase.}
% The above constraint can be represented as a linear constraint in the following way: $r_{k,t}-r_{k,t-1}=\text{diff}^{+}_{k,t}-\text{diff}^{-}_{k,t}$ where $\text{diff}_{k,t}^{+}$ ($\text{diff}_{k,t}^{-}$, resp.) is the positive part (negative part, resp.) of the difference $r_{k,t}-r_{k,t-1}$. Thus, the above constraint becomes
% \begin{eqnarray}
% \sum_{t=t_k}^{t_{dead}-1}\text{diff}^{+}_{k,t}+\text{diff}^{-}_{k,t}-l\leq BU, \quad \sum_{t=t_k}^{w_i-t_k-1}\text{diff}^{+}_{i,t}-\text{diff}^{-}_{i,t}-N_i\leq BU_i \forall i\in \mathcal{K}_0.\label{eq:battery_utilization}\\
% r_{k,t}-r_{k,t-1}=\text{diff}^{+}_{k,t}-\text{diff}^{-}_{k,t}, \quad r_{i,t}-r_{i,t-1}=\text{diff}^{+}_{i,t}-\text{diff}^{-}_{i,t-1} \quad \forall i \in \mathcal{K}_0.\label{eq:difference}
% \end{eqnarray}
 \section{Problem Formulation}
  \subsection{User's utilities} \label{sec:userutility}User's utility for the contract $(l,t,BU)$ is denoted as $U^{BU}_{k,l,t}$ which is a random variable. The realized value  $u^{BU}_{k,l,t}$ is only known to the user $k$, but not known to the charging station in general.  The payoff of user $k$ or {\em user's surplus} if she selects the  contract $(l,t,BU)$ is $u^{BU}_{k,l,t}-p^{BU}_{k,l,t}$ (Fig.~\ref{fig:ev_charging}). If she rejects all her payoff is $0$. The user will select the contract that fetches the maximum payoff to her. 
  
 Thus, for a menu of  prices $p^{BU}_{k,l,t}$, the user $k$ selects $A^{BU}_{k,l,t}\in [0,1]$ such that it maximizes the following
\begin{align}
& \text{maximize} \sum_{l=1}^{L}\sum_{t=t_k+1}^{T}\sum_{BU=0}^{BU_{max}}A^{BU}_{k,l,t}( u_{k,l,t,BU}-p^{BU}_{k,l,t})\nonumber\\
& \text{subject to } \sum_{l=1}^{L}\sum_{t=t_k+1}^{T}\sum_{BU=0}^{BU_{max}}A_{k,l,t,BU}\leq 1
\end{align}
 Note that $A^{BU}_{k,l,t}>0$ only if $\max_{i,j,b}\{u_{k,i,j,b}-p_{k,i,j,b}\}=u^{BU}_{k,l,t}-p^{BU}_{k,l,t}$.  If such a solution is not unique, any convex combination of these solutions is also optimal since a user can select any of the maximum payoff contracts. 
 
 We denote the  decision as $A^{BU}_{k,l,t}(\mathbf{p}_{k})$. Note that the decision whether to accept the menu price $p^{BU}_{k,l,t}$  depends not only on the price $p^{BU}_{k,l,t}$ but also other price menus i.e. $p^{b}_{k,i,j}$ where $i\in \{1,\ldots, L\}$ and $j\in \{t_k+1,\ldots, T\}$  and $b\in \{0,\ldots,BU_{max}\}$. This is because the user only selects the price menu which is the  most favorable. Note that if the maximum payoff that user gets among all the price menus (or, contracts) is negative, then the user will not charge i.e., $A^{b}_{k,l,t}=0$ for all $l$, $t$ and $b$. We also assume that if there is a tie between charging and not charging, then the user will decide to charge i.e., if the maximum payoff that user can get is $0$, then the user will decide to charge \footnote{However, our result can be readily extended to the other options, in that case the price strategies given in this paper have to decreased by an $\epsilon>0$ amount.}.
 
  %We also assume that $C_k(\cdot)$ is random whose value is again from some distribution.  

 %\textbf{Utilization of the batteries of EVs} Frequent charging and discharging may degrade the battery life of the EVs. The users have to be paid adequately for such undesirable events. Hence, the total cost of the charging station is
% \begin{align}
% \mathcal{P}_{l,t_{dead}}:& \text{minimize } \sum_{t=t_{k}}^{T-1}c_tq_t-g_tx_t+\alpha(|r_{k,t}-r_{k,t-1}|-l)+\sum_{i\in\mathcal{K}_0}\alpha(|r_{i,t}-r_{i,t-1}|-N_i)\nonumber\\
% & \text{subject to } (\ref{eq:charged})-(\ref{eq:storage})\nonumber\\
% & \text{var } q_t\geq 0, x_t\geq 0, r_{k,t}, r_{i,t}\forall i\in \mathcal{K}_0.
% \end{align}
% For user $k$, $r_{k,t_k-1}$ is assumed to be $0$. Note that an user is only compensated when the battery is used for both battery charging and discharging during the course of the stay. If it consists of only type of cycle (either charging or discharging) the user will not be compensated. $\alpha$ is computed depending on the total charging and discharging cycles the battery can handle. 

\vspace{-0.1In}
\subsection{Myopic Charging Station}
 Since the users arrive to the charging station at any time throughout the day, the charging station does not know the exact arrival times for the future vehicles.  We consider that the charging station is myopic or near-sighted  i.e., it selects its price for user $k$ without considering the future arrival process of the vehicles. However, it will consider the cost incurred to charge the existing EVs. Note that as the number of existing users increases, the marginal cost can increase to fulfill a contract for an arriving user, hence, such a pricing strategy may not maximize the payoff in a long run. Note that, {\em a myopic pricing strategy is optimal in the case the marginal cost of fulfilling a demand of a new user  is independent of the number of existing users.}  

 In practice, the charging station often has fixed number of charging spots. Thus, the charging station may want to select high prices for user $k$, in order to make the charging spots available for the users who can pay more but only will arrive in future\footnote{The above consideration is left for the future work.}. Hence, considering the future arrival process of the vehicles may increase the profit of the charging station. However, such a pricing strategy is against the first come first serve basis which is the current norm for charging vehicles. Our approach can be seen as a fair allocation process, where the charging station serves users based on the first come first serve basis. % Later in Section~\ref{sec:simulation_results}, we show that since the charging station can control the time spent by a EV through pricing, our approach requires lower number of charging spots compared to the existing pricing mechanism.
 
  \subsection{Optimal Cost of the charging station to fulfill a menu}  The optimal cost for the charging station to fulfill the contract $(l,t_{dead},BU)$ to the user $k$ is thus
 \begin{align}\label{eq:vlt}
 \mathcal{P}^{BU}_{l,t_{\text{dead}}}:& \text{minimize } \sum_{t=t_k}^{T-1}(c_tq_t-g_tx_t)\nonumber\\
 & \text{subject to } (\ref{eq:charged}), (\ref{eq:charge}),(\ref{eq:chargeanddischarge}),(\ref{eq:chargeeachev}), (\ref{eq:dischargelimit}), (\ref{eq:limit}), (\ref{eq:battery_level}), (\ref{eq:battery_capacity}), (\ref{eq:storage}), (\ref{eq:battery_utilization}).\nonumber\\
 & \text{var } q_t\geq 0,  x_t\geq 0, e_t\geq 0, s_t\geq 0, \nonumber\\
 & r^{+}_{i,t}\geq 0, r^{-}_{i,t}, r_{i,t}\geq0 \quad \forall i\in \mathcal{K}_0\cup \{k\}.
 \end{align}
 where $c_t$ is the per unit cost of buying energy from the grid and $g_t$ is the per unit cost of selling energy to the grid. Note that our model can also incorporate time varying, strictly increasing convex costs $C_t(\cdot)$ or concave prices $G_t(\cdot)$.We assume that $c_t\geq g_t$ in order to avoid arbitrage opportunity. However, the above cost structures will destroy the properties of the linear programming. 
 
 Our model can be easily extended to the setting where there are upper bounds on $q_t$ and $x_t$. The problem will still be a linear programming problem.
 
 The cost $c_t$ and the price $g_t$ are assumed to be known. In case they change in a dynamic manner i.e., they are real time prices, we consider that $c_t$ is the expected cost, and $g_t$ is the expected price. 
 
 Note that the above problem is a linear programing problem (similar to the setting considered in \cite{tech_ev}) and thus, it is easy to solve. However, the number of constraints are hugher compared to \cite{tech_ev} as we consider the V2G service. 
 
 In the following theorem, we show that the optimal cost  to the charging station given in (12), $r_{i,t}^+r_{i,t}^- =0$ and $e_ts_t=0$ for all $t$, and $i\in\mathcal{K}_0\cup \{k\}$.
 \begin{theorem}\label{thm:noposinegative}
 In an optimal solution of the problem $\mathcal{P}^{BU}_{l,t_{\text{dead}}}$, we must have the decision variables $r^{+}_{i,t}r^{-}_{i,t}=0$ $\forall i\in \mathcal{K}_0\cup \{k\}$  and $t$. Further, in an optimal solution, we also have $e_ts_t=0$ for all $t$.
 \end{theorem}
 {\em Proof}: See Appendix~\ref{sec:proof_thmnoposinegative}.\qed
 
 The above theorem shows that in an optimal solution, neither the  EV nor the battery of the charging station simultaneously charge and discharge. 

Now we calculate the {\em additional} cost imposed to the charging station for fulfilling the new contract. First, we introduce some notations which we use throughout. 
\begin{definition}\label{defn:vlt}
The charging station has to  incur the cost $v^{BU}_{l,t_{\text{dead}}}$ for serving the existing users and the contract $(l,t_{\text{dead}},BU)$ for the new user $k$ where $v^{BU}_{l,t_{\text{dead}}}$ is the value of the linear optimization problem $\mathcal{P}^{BU}_{l,t_{\text{dead}}}$.
\end{definition}
 Since $\mathcal{P}^{BU}_{l,t_{\text{dead}}}$ is a linear optimization problem, it is easy to compute $v^{BU}_{l,t}$. Further, note that if the above problem is infeasible for some $l$, $t$ and $BU$, then we consider $v^{BU}_{l,t}$ as $\infty$. 

\begin{definition}\label{defn:v-k}
Let $v_{-k}$ be the amount that the charging station has to incur to satisfy the requirements of the existing EVs if the new user does not opt for any of the price menus. 
\end{definition}
If user $k$ does not accept any price menu, then the charging station still needs to satisfy the demand of existing users i.e., the charging station must solve the problem $\mathcal{P}^{BU}_{l,t_{\text{dead}}}$ with $r_{k,t}=0$. $v_{-k}$ is the value of that optimization problem.   From Definitions~\ref{defn:vlt} and \ref{defn:v-k} we can visualize $v^{BU}_{l,t_{\text{dead}}}-v_{-k}$ as the additional cost or marginal cost to the charging station when the user $k$ accepts the price menu $p^{BU}_{k,l,t_{\text{dead}}}$.  

We assume that the prediction $\bar{E}^{t}$ is perfect for all future times and is known to the charging station. However, menu-based pricing approach can be extended to the setting where the estimated generation does not match the exact amount. First, we can consider a conservative approach where $\bar{E}^t$ can be treated as the worst possible renewable energy generation. As a second approach, we can  accumulate various possible  scenarios of the renewable energy generations, and try to find the cost to fulfill a contract for each such scenario. For example, if there are $M$ number of possible instances of the renewable energy generation amount in future. Then, we can find the optimal cost for each such instance of renewable energy generation $\bar{E}^{m,t}$ where $m\in \{1,\ldots, M\}$ instead of $\bar{E}^t$. We then can compute the average (or, the weighted average, if some instance has greater probability) of the optimal costs, and that cost can be taken as the cost of fulfilling a certain contract.

The following result entails that the V2G service indeed reduces the cost of the charging station. 
\begin{lemma}\label{lm:v2g}
$v^{BU_1}_{l,t}\leq v^{BU_2}_{l,t}$ for any $BU_1>BU_2$. Thus, $v^{BU}_{l,t}\leq v^{0}_{l,t}$ for any $BU\geq 1$.
\end{lemma}
{\em Outline of the Proof}: Note from (\ref{eq:battery_utilization}) that the decision space is more restricted for $BU_2$ compared to $BU_1$ if $BU_1>BU_2$. Every feasible solution with $BU_2$ is also feasible with $BU_1$ since $BU_1>BU_2$. Thus, the above result follows.\qed

Since the charging station can sell some energies taken from the EVs, the cost of fulfilling the contract $(l,t,BU)$ will be lower as $BU$ will increase. 
\subsection{Charging Station's Profit}\label{sec:profit}
If the user $k$ selects the menu $(l,t,BU)$, then the additional cost incurred by the charging station is $v^{BU}_{l,t}-v_{-k}$. Hence, the profit of the charging station is $p^{BU}_{k,l,t}-v^{BU}_{l,t}+v_{-k}$. If the user does not select any of the menus, then the charging station gets $0$ profit. Thus, the expected profit of the charging station is \begin{align}\label{eq:profitmax_expec}
\sum_{BU=0}^{BU_{max}}\sum_{l=1}^{L}\sum_{t=t_k+1}^{t=t_k+T} (p^{BU}_{k,l,t}-v_{l,t,BU})\Pr(R^{BU}_{k,l,t})
\end{align}
where $\Pr(R^{BU}_{k,l,t})$ is the probability of the event that the menu $(l,t,BU)$ is accepted by the user $k$. Note from (\ref{eq:profitmax_expec})  that the profit maximization is a difficult problem as the menu selected by an user inherently depends on the prices selected for other menus. For example, if the price selected for a particular contract is high, the user will be reluctant to take that as compared to a lower price one. The profit is a discontinuous function of the prices and thus, the problem may not be convex even when the marginal distribution of the utilities are concave.

\subsection{Objectives}
We consider that the charging station decides the price menus in order to fulfill one of the two objectives (or, both)-- i) Social Welfare Maximization and ii) its profit maximization.
\subsubsection{Social Welfare}
The social welfare is the sum of  user surplus and the profit of the charging station. As discussed in Section~\ref{sec:userutility} for a certain realized values $u^{BU}_{k,l,t}$ if the user $k$ selects the price menu $p^{BU}_{k,l,t}$, then its surplus is $u^{BU}_{k,l,t}-p^{BU}_{k,l,t}$, otherwise it is $0$.

As discussed in Section~\ref{sec:profit} the profit of the charging station is $p^{BU}_{k,l,t}-v^{BU}_{l,t}+v_k$ for a given price $p^{BU}_{k,l,t}$  if the user selects the menu, and it is $0$ otherwise.  Hence, the social welfare maximization problem is to select the price menu $p_{k,l,t}$ which will maximize the following
\begin{align}\label{eq:socialwelfare_maxproblem}
 \mathcal{P}_{\text{perfect}}: & \text{max} \sum_{l=1}^{L}\sum_{t=t_k+1}^{T}\sum_{BU=0}^{BU_{max}}(u^{BU}_{k,l,t}-v^{BU}_{l,t}+v_{-k})A_{k,l,t}(\mathbf{p}_k)\nonumber\\
& \text{var }: p^{BU}_{k,l,t}\geq 0.
\end{align}
Recall that in order to find $v^{BU}_{l,t}$ we have to solve $\mathcal{P}_{l,t}$(cf.~(\ref{eq:vlt})) which is a constrained optimization problem. 
%Since EV is expected to increase the social value such as providing cleaner environment, and higher energy efficiency, hence, it is  important for regulator (e.g. FERC) whether there exists a pricing strategy which maximizes the social welfare of the system. Being a social planner, the FERC may want the system becomes socially efficient. If the charging station is operated by the regulator or some government agency, then the main objective is indeed maximizing the social welfare or user surplus is maintained.

Since the charging station is unaware of the utilities of the users, the charging station has two choices-i) decides a price and hopes that it will maximize the social welfare for the realized values of utilities ({\em ex-post} maximization) , or ii) decides a price and hopes that it will maximize  the social welfare in an expected sense ({\em ex-ante} maximization). % then the charging station may want to maximize either in  i) ex-ante  or ii) ex-post  social welfare.
%\begin{definition}\label{defn:ex-ante}
%In the {\em ex-ante social welfare} maximization, the social welfare is maximized in the expected sense where the expectation is taken over the joint distributions of the random variables $U_{k,l,t}$.
%\end{definition}
%\begin{definition}\label{defn:ex-post}
%In the {\em ex-post social welfare} maximization, the social welfare is maximized for each possible realization of the random variables $U_{k,l,t}$. 
%\end{definition}
Thus,  the {\em ex-ante} maximization  does not guarantee that the social welfare will be maximized for every realization of the random variables $U_{k,l,t}$. However, in the {\em ex-post} maximization, the social welfare is maximized for each possible realization of the random variables. Thus, {\em ex-post} maximization is a stronger concept of maximization (and thus, also more desirable) and it is not necessary that there exist pricing strategies which maximize the ex-post  social welfare. However, we show that in our setting there exist pricing strategies which maximize the ex-post social welfare. Note that {\em ex-post} social welfare maximization is the same as (\ref{eq:socialwelfare_maxproblem}).

\subsubsection{Profit Maximization}
Social welfare maximization does not guarantee that the charging station may get a positive profit. It is important for the wide scale deployment of the charging stations that the charging station must have some profit. The charging station needs to select $p^{BU}_{k,l,t}$ in order to maximize the expected profit (given in (\ref{eq:profitmax_expec})).

Note that in order to select optimal $p^{BU}_{k,l,t}$, the charging station has to obtain $v^{BU}_{l,t}$ i.e., it has to solve the problem $\mathcal{P}^{BU}_{l,t}$ (cf.~(\ref{eq:vlt})) for each choice of $l$, $t$ and $BU$.

\subsubsection{Separation Problem}
Note that in order to select optimal $p^{BU}_{k,l,t}$, the charging station has to obtain $v^{BU}_{l,t}$  and $v_{-k}$ (Definitions~\ref{defn:vlt} \& \ref{defn:v-k}). However,  $v_{l,t}$ and $v_{-k}$  do not depend on $p^{BU}_{k,l,t}$. Hence, we can separate the problem--first the charging station finds $v^{BU}_{l,t}$ and $v_k$, and then it will select $p^{BU}_{k,l,t}$ to fulfill the objective. {\em We now focus on finding optimal} $p^{BU}_{k,l,t}$.

 \section{Result: Ex-post Social Welfare Maximization}
First, we state the optimal values of the social welfare for any given realization of the user's utilites. Next, we state a pricing strategy which attains the above optimal value. 

Note that if $u^{BU}_{k,l,t}-v^{BU}_{l,t}+v_{-k}<0$ for each $l$,$t$, and $BU$,  then the social welfare is maximized when the user $k$ does not charge. In this case, the optimal value of social welfare is $0$.

On the other hand if $u^{BU}_{k,l,t}-v^{BU}_{l,t}\geq -v_{-k}$ for  some $l$, $t$, and $BU$ then the  social welfare is maximized when  the user $k$ charges its car. If the user accepts the price menu $p^{BU}_{k,l,t}$, then the social welfare is $u^{BU}_{k,l,t}-v^{BU}_{l,t}+v_{-k}$. Thus, the maximum social welfare in the above scenario is $\max_{l,t,BU}(u^{BU}_{k,l,t}-v^{BU}_{l,t}+v_{-k})$ as described in the following theorem--
\begin{theorem}\label{thm:max_socialwelfare}
The maximum value of social welfare is $\max\{\max_{l,t,BU}(u^{BU}_{k,l,t}-v^{BU}_{l,t}+v_{-k}),0\}$.
\end{theorem}
 
 %The ex-post social welfare is the sum of the realized value of the utility of the user and the profit of the charging station. 
 %Thus, for realized values of the utilities of user $k$, the ex-post social welfare is 
% \begin{align}
% \sum_{BU=0}^{BU_{max}}\sum_{l=1}^{L}\sum_{t=t_k+1}^{T}(u^{BU}_{k,l,t}-v_{l,t,BU})\mathbbm{1}(R^{BU}_{k,l,t})
% \end{align}
 We obtain
 \begin{theorem}\label{thm:max_socv2g}
 If $p^{BU}_{k,l,t}=v^{BU}_{l,t}-v_k$ , then it will maximize the ex-post social welfare.
 \end{theorem}
 {\em Outline of the Proof}: If the price is set according to the above, the user will not select any contract if $\max_{l,t,BU}(u^{BU}_{k,l,t}-v^{BU}_{l,t}+v_{-k})<0$ which gives $0$ social welfare. The user select the contract which maximizes the payoff if $\max_{l,t,BU}(u^{BU}_{k,l,t}-v^{BU}_{l,t}+v_{-k})\geq 0$. In this case, the ex-post social welfare is $\max_{l,t,BU}(u^{BU}_{k,l,t}-v^{BU}_{l,t}+v_{-k})$. Hence, the result follows.\qed
 
 However, the above pricing mechanism does not give any non-zero profit to the charging station. The above pricing strategy is distribution independent, it holds for {\em any arbitrary distribution. }
 
 Also note that the pricing strategy also maximizes the social welfare in the long run when the additional cost of fulfilling a contract (i.e. $v^{BU}_{l,t}-v_{-k}$) does not depend on the existing users in the charging station. 
%Hence,
%\begin{corollary}\label{cor:longrun}
%The pricing strategy $p_{k,l,t}=v_{l,t}-v_{-k}$ maximizes the ex-post social welfare over any time interval $T$ even when the charging station is not myopic if $v_{l,t}-v_{-k}$ does not depend on the existing users for all $l$ and $t$
%\end{corollary}
The condition that $v_{l,t}-v_{-k}$ is independent of the existing EVs in the charging station is satisfied if either all demand can be fulfilled using renewable energy or there is no renewable energy generation. Hence, in the {\em two above extreme cases, the myopic pricing strategy is also optimal in the long run}.
 \section{Result:Profit Maximization}
 \subsection{Maximum Profit under ex-post social welfare maximization}
 %Note from Theorem~\ref{thm:profit_max} that the profit maximization pricing strategy which maximizes the social welfare requires that the charging station has the complete information of the utilities of the users which is not possible in this case. Hence, such a pricing strategy can not be implemented when the charging station does not know the exact utilities of the users. 

 %
%First, we provide the expression of the expected profit of the charging station

We have already seen  a pricing strategy which maximizes the ex-post social welfare in Theorem~\ref{thm:max_socv2g}, however, this pricing strategy does not give any positive profit. Naturally the question arises what is the pricing strategy that maximizes the expected profit of the charging station  which will also maximize the ex-post social welfare.

 %The utilities $U_{k,l,t}$ are also dependent as $U_{k,l+1,t}\geq U_{k,l,t}$ since higher energy within the same deadline should induce higher utility, hence, the standard techniques to finding prices for independently distributed random variables can also not be applied. 

We show that there exists a pricing strategy which may provide better profit to the charging station while maximizing the ex-post social welfare. First, we introduce a notation which we use throughout.
\begin{definition}\label{defn:lower_endpoint}
Let $L^{BU}_{k,l,t}$ be the lowest end-point of the marginal distribution of the utility $U^{BU}_{k,l,t}$.
\end{definition}
\begin{theorem}\label{thm:profitmax_uncertainty}
Consider the pricing strategy: 
\begin{align}\label{eq:pricesocialmax}
p^{BU}_{k,l,t}=v^{BU}_{l,t}-v_{-k}+(\max_{i,j,b}\{L^{b}_{k,i,j}-v^{b}_{i,j}+v_{-k}\})^{+}.
\end{align} The pricing strategy maximizes the ex-post social welfare.

The profit  is $(\max_{i,j,b}\{L^{b}_{k,i,j}-v^{b}_{i,j}+v_{-k}\})^{+}$.
\end{theorem}
{\em Outline of proof}: First, note that adding a constant does not change the optimal solution. Hence, if  $(l^{*},t^{*},b^{*})=\text{argmax}_{l,t,b}(u^{b}_{k,l,t}-v^{b}_{l,t}+v_{-k})$, then $(l,^{*},t^{*}) $ is also optimal for price strategy in (\ref{eq:pricesocialmax}). Now, if the price is set according to (\ref{eq:pricesocialmax}),  when $u^{b^{*}}_{l^{*},t^{*}}\geq v^{b^{*}}_{l^{*},t^{*}}-v_{-k}$, the user always select the contract. If the condition is not satisfied, the user does not select the contract. Hence, the ex-post social welfare is always maximized. \qed %The  rest of the proof follows from Lemma~\ref{lm:delta}. Lemma~\ref{lm:delta} entails that there exists some $\delta>0$ such that it will ensure that the price strategy is off from the social welfare maximizer pricing strategy by at most $1-\epsilon$ in probability.

%Also note  the similarity with Theorem~\ref{thm:profit_max}. If the user knows the utility, then $L_{k,l,t}=u_{k,l,t}$ as there is no uncertainty and we get back the pricing strategy stated in Theorem~\ref{thm:profit_max}. 

Note that if $\max_{l,t,BU}(L^{BU}_{k,l,t}-v^{BU}_{l,t}+v_{-k})>0$, then such a pricing strategy will provide a positive profit to the charging station. Hence, if $v^{BU}_{l,t}-v_{-k}$ is smaller than the lowest end-point of the distribution function, then the profit will be positive. % Note that if the harvested energy is large, then the V2G service from the EVs may not be required, as the charging station may have enough energy to sell to the grid and serve the users.  If the harvested energy is not very high, then the charging station may get additional profit if the users participate in the V2G service. Hence, {\em when the harvested energy is not enough, the charging station can get additional profit from the V2G service.}%Thus, the above illustrates the importance of the storage and harvesting energy devices in the charging station. The regulator (e.g. FERC) can also provide incentives to the charging station to set up those devices as the pricing strategy can give a positive profit to the charging station as well as maximize the ex-post social welfare. 

%In the extreme, when $v_{l,t,BU}=0$ for all $l$ and $t$, then the maximum profit that the charging station makes is $\max_{l,t,BU}\{L^{BU}_{k,l,t}\}$ while maximizing the ex-post social welfare. It also increases the user  surplus, as the price set by the charging station decreases. Thus, the impact of higher degrees of renewable energy  integration for the charging station increases both the profit of the charging station and the user surplus. However,  further decreasing $v_{l,t,BU}$ will not have any effect on the profit of the charging station as well as the user surplus, thus, it also shows the investment that the charging station needs to make for storage and renewable energy harvesting devices.
% Thus,
%\begin{corollary}
%Both the user's surplus and the profit of the charging station increase with increase in the renewable energy generation.
%\end{corollary}

Also note that the users which have higher utilities i.e., higher $L^{BU}_{k,l,t}$ will give more profits to the charging station. 

The charging station needs to know the lowest end-points of the support set of the utilities unlike in Theorem~\ref{thm:max_socv2g}. However, the charging station does not need to know the exact distribution functions of the utilities similar to Theorem~\ref{thm:max_socv2g}. The lowest end-point can be easily obtained from the historical data. 

The pricing strategy maximizes the ex-post social welfare similar to Theorem~\ref{thm:max_socv2g}. This is also the {\em maximum possible profit that the charging station can have under the condition that it  maximizes the ex-post social welfare with probability $1$}. However, it may not maximize the expected profit of the charging station. In other words, the pricing strategy which maximizes the expected profit needs not maximize the ex-post social welfare.  %Hence unlike the scenario where the charging station is clairvoyant here there may not exist a profit maximization strategy which is also a social welfare maximizer. Note that the user surplus is not $0$, hence {\em uncertainty regarding the values of the utilities  help the users to gain positive surpluses} unlike in the scenario where the charging station is clairvoyant.

\textbf{When the charging station is clairvoyant}: We have seen that if the charging station is unaware of the realized values of the utilities, there is no pricing strategy which maximizes both the expected profit and the ex-post social welfare. However, we show that  if the charging station is clairvoyant {\em i.e.}, it is aware of the  realized values of the utilities of the users, then there exists a pricing strategy which both maximizes the social welfare and the profit of the charging station. 

First, we introduce a notation.
\begin{definition}\label{defn:optimal}
Let $(l^{*},t^{*},b^{*})=\arg\max_{l,t,b}\{u^{b}_{k,l,t}-v_{l,t}\}$.
\end{definition}
\begin{lemma}\label{thm:profit_max}
Let $p^{BU}_{k,l,t}=v^{BU}_{l,t}-v_{-k}+(u^{b^{*}}_{k,l^{*},t^{*}}-v^{b^{*}}_{l^{*},t^{*}}+v_{-k})^{+}$ where $(l^{*},t^{*},b^{*})$ is given in Definition~\ref{defn:optimal}. Such a pricing strategy maximizes the profit as well as the social welfare.
\end{lemma}
{\em There can be other  pricing strategies which simultaneously maximize the social welfare and the profit.} %For example, if $p_{k,l,t}$ is $\infty$ for all $(l,t)\neq (l^{*},t^{*},b^{*})$ and $p^{b^{*}}_{k,l^{*},t^{*}}=v^{b^{*}}_{l^{*},t^{*}}-v_{-k}+(u^{b^{*}}_{k,l^{*},t^{*}}-v^{b^{*}}_{l^{*},t^{*}}+v_{-k})^{+}$, then it also maximizes the profit of the charging station. Thus, when the charging station is clairvoyant, it can give only one possible contract to the EVs.
Though the joint  profit and social welfare maximizing pricing strategy may not be unique, the profit of the charging station is the unique and  is given by 
\begin{align}
\max\{u^{b^{*}}_{k,l^{*},t^{*}}-v_{l^{*},t^{*}}+v_{-k},0\}
\end{align}
The above pricing strategy is an example of {\em value-based} pricing strategy where prices are set depending on the valuation or the utility of the users \cite{Hinterhuber}. In contrast, the price strategy stated in Theorem~\ref{thm:max_socv2g} is an example of {\em cost-based} pricing strategy where the prices only depend on the costs. In the value-based pricing strategy, the user surplus decreases, in fact it is\footnote{If the user is reluctant to charge if it does not get a positive payoff, then, we can reduce the price by $\epsilon>0$ amount. In that case, it will be $(1-\epsilon)$ optimal profit maximizing strategy.} $0$ in our case.  Thus all the user surplus is transferred as the profit of the charging station. {\em Thus, uncertainty regarding the utilities enhance the user's surplus.}

\subsection{Guaranteed positive profit to the Charging station}\label{sec:price_uncertainty}
Theorem~\ref{thm:profitmax_uncertainty} entails that the charging station only has a positive profit if $\max_{l,t,b}\{L^{b}_{k,l,t}-v^{b}_{l,t}+v_k\}>0$. If the above condition is not satisfied, the charging station's profit will be $0$. In the following we consider a pricing strategy which will give a guaranteed positive profit to the charging station. 

%First note that by the continuity of the joint distribution function we have the following
%\begin{lemma}\label{lm:delta}
%Let for each $\epsilon>0$, there exists a $\delta>0$ such that 
%\begin{align}
%& \Pr(\max_{l,t,BU}\{U^{BU}_{k,l,t}-v_{l,t}+v_{-k}\}\geq 0)\nonumber\\
%& \leq \epsilon+\Pr(\max_{l,t,BU}\{U^{BU}_{k,l,t}-v_{l,t}+v_{-k}-\delta\}\geq 0)
%\end{align} 
%\end{lemma}
%
%\begin{theorem}\label{thm:approx}
%Fix an $\epsilon>0$. Now, consider the pricing strategy
%\begin{align}\label{eq:approx_price}
%p^{BU}_{k,l,t}=v_{l,t,BU}-v_{-k}+\delta(\epsilon)
%\end{align}
%where $\delta(\epsilon)$ is the $\delta$ which satisfies the Lemma~\ref{lm:delta}.
%
%Then such a pricing strategy maximizes  the ex-post social welfare with probability $1-\epsilon$. 
%\end{theorem}
Consider the pricing strategy
\begin{align}\label{eq:guranteed_profit}
p^{BU}_{k,l,t}=v^{BU}_{l,t}-v_{-k}+\beta%+(\max_{i,j,b}\{L^{b}_{k,i,j}-v^{b}_{i,j}+v_{-k}\})^{+}
\end{align}
where $\beta>0$. %Note that  this pricing strategy can be readily implemented.

Note that the pricing strategy stated in (\ref{eq:pricesocialmax}) is a variant of the pricing strategy stated in (\ref{eq:guranteed_profit}) where $\beta=(\max_{i,j,b}\{L^{b}_{k,i,j}-v^{b}_{i,j}+v_{-k}\})^{+}$ if we allow $\beta$ can also be $0$.

The pricing strategy stated in (\ref{eq:guranteed_profit}) gives the same positive profit irrespective of the menu selected by the user. The regulator such as FERC can select a $\beta$ judiciously to trade off between the profit of the charging station and the social welfare. From Lemma~\ref{lm:v2g}, the pricing strategy stated in  (\ref{eq:guranteed_profit}) selects {\em lower price for higher battery utilization.} This is desirable, as the user needs to be given incentive to battery utilization.  

Higher $\beta$ will deter the user's surplus. Also note that when $\beta>0$, it may not maximize the ex-post social welfare from Theorem~\ref{thm:profitmax_uncertainty}. Very high value of $\beta$ also decreases the profit of the charging station, as users will be reluctant to accept any of the menus. %Note that in the above pricing strategy, irrespective of the menu selected, the charging station always gets a uniform $\delta(\epsilon)$ amount of profit if the user selects a menu. 

%   Note also that since it also gives a positive amount of profit to the charging station and maximizes the social welfare with a high probability, the above pricing strategy is also preferable for both the regulator and charging station operator.
%
%Note that {\em the assumption of continuous distribution is key}. If the distributions are not discrete, then $\delta(\epsilon)$ may be $0$. Hence, the charging station may  get zero profit even in this case.

The expected profit of the charging station for the above pricing strategy is--
\begin{theorem}\label{thm:expected_payoff}
The expected profit of the charging station when it selects price according to (\ref{eq:guranteed_profit}) is $\beta\max_{l,t,BU}\{\Pr(U^{BU}_{k,l,t}\geq v^{BU}_{l,t}-v_{-k}+\beta)\}$.
\end{theorem}
{\em Outline of the Proof}: Note that if a user selects any of the contracts, then the charging station's profit is $\beta$. Hence, the charging station's expected profit is $\beta$ times the probability that at least one of the contracts will be accepted.\qed

Now, we provide an example where the pricing strategy in (\ref{eq:guranteed_profit}) can also maximize the expected profit for a suitable choice of $\beta$.  First, we introduce a notation
\begin{definition}\label{defn:alpha}
Let $\zeta=\max\{\gamma| \gamma\in \text{argmax}_{\beta\geq 0}\beta\{\max_{i,j,b}\Pr(U^{b}_{k,i,j}\geq \beta+v^{b}_{i,j}-v_{-k}\}\}$. 
\end{definition}
Note that since $U^{BU}_{k,l,t}$ is bounded and the probability distribution is continuous, thus, $\zeta$ exists. Note from Theorem~\ref{thm:expected_payoff} that $\zeta$ corresponds to $\beta$  the charging station can get the maximum possible expected profit when the prices are of the form (\ref{eq:guranteed_profit}).  If the utilities are drawn from a strictly increasing continuous distribution, then the set of $\gamma$ would be singleton and we do not need to specify the maximum. 

Now, consider the pricing strategy 
\begin{align}\label{eq:price_alpha}
p^{BU}_{k,l,t}=v^{BU}_{l,t}-v_{-k}+\zeta.
\end{align}
where $\zeta$ is as given in Definition~\ref{defn:alpha}.
The above pricing strategy maximizes the profit for a class of utility functions which we describe below.
\begin{assumption}\label{assum:utility}
Suppose that the utility function $U^{BU}_{k,l,t}=(Y^{BU}_{k,l,t}+X_k)$ for all $l$, $t$ \& $BU$; $Y^{BU}_{k,l,t}$ is a constant and known to the charging station, and $X_k$ is a random variable whose realized value is not known to the charging station.
\end{assumption}
In the above class of utility function, the uncertainty is only regarding the realized value of the random variable $X_k$. Note that $X_k$ is independent of $l$, $BU$ and $t$, hence,$X_k$ is considered to be an additive noise.

It is important to note that we do not put any assumption whether {\em $X_k$ should be drawn from a continuous or discrete distribution.} However, if the distribution is discrete, we need the condition that $\zeta$ must exist. 
%The charging station is fully aware of the $Y_{k,l,t}$. In the above class of utilities, if $p_{k,l_1,t_1}-p_{k,l_2,t_2}>Y_{k,l_1,t_1}-Y_{k,l_2,t_2}$, then the charging station will never select the price menu $p_{k,l_1,t_1}$. Hence, in the above class of utility functions, the charging station is aware of the price differences that it should select so that the user selects that option despite the uncertainty due to $X_k$. 

\begin{theorem}\label{thm:aclassutility}
The pricing strategy stated in (\ref{eq:price_alpha}) maximizes the expected profit of the charging station (given in (\ref{eq:profitmax_expec})) when the utility functions are of the form given in Assumption~\ref{assum:utility}. 
\end{theorem}
%{\em Proof}: See Appendix~\ref{sec:proof_thmaclassutility}.\qed

{\em The above result is surprising. It shows that a simple pricing mechanism such as the fixed profit can maximize the expected payoff for a large class of utility functions.}
  However, if the utilities do not satisfy Assumption~\ref{assum:utility} the above pricing strategy may not be optimal. 

\section{The user's participation in the V2G service and the profitability}\label{sec:v2g_profitability}
The V2G service will proliferate only if the users participate in that service. The charging station can attain extra profits through the V2G service. However,  the users will only select the menu with positive battery utilization if they get enough compensation. Thus, the charging station's profit inherently depends on whether the users have incentives to participate in the V2G services.  In this section, we will analyze the conditions under which the users will be willing to participate in the V2G service, and the profit of the charging station will increase. 
\subsection{Cost of Battery Utilization}
First, we discuss the cost of battery utilization. Users will strictly prefer lower utilization as lower $BU$ will increase the battery life. A higher battery utilization may increase the battery degradation cost \cite{Delucchi,zhou}. We denote the cost associated with the utilization $BU$  for user $k$ is $C_k(BU)$ where 
$C_k(\cdot)$ is a strictly increasing function. 
 
 The cost $C_k(\cdot)$ depends on the the state of the battery\footnote{The cost may also depend on the total number of charging and discharging cycles the car has gone through. It may also depend on the user's willingness to participate in the V2G service.} \cite{Delucchi,zhou}.  We assume that the user's utility $U^{BU}_{k,l,t}$ is
 \begin{align}
 U^{BU}_{k,l,t}=U_{k,l,t}-C_k(BU)
 \end{align}
We consider that the cost function $C_k(BU)$ for the battery utilization as a linear function i.e. $C_k(BU)=\alpha_k BU$.  Recently, \cite{zhou} shows that the per unit degradation cost for discharging remains almost constant for a wide range of values. Hence,  a linear cost model can be a good approximation of the cost function. However, our analysis can be easily extended to other cost models. The charging station and even the user may not know the exact value of $\alpha_k$. But, the EV manufacturer can easily provide the {\em pessimistic} approximation of $\alpha_k$ such as the worst possible battery degradation cost for per unit of energy.\footnote{Recently, \cite{zhou} shows that the battery degradation cost of Li-Ion battery for per unit of energy is shown to be between $4$ cents and $7$ cents. In this example $\alpha_k$ can be taken as $7$ cents per kwh. } 

The realized value of the utility function of user $k$  is now $u_{k,l,t}-\alpha_k BU$.  %Note that the rest of the analysis can also be extended to other cost functions. %We assume that $U_{k,l,t}$  is a random variable whose realized value $u_{k,l,t}$  is only known to the user. We assume that $U_{k,l,t}$ is drawn from a joint distribution over $l$ and $t$ which may not be independent. $C_k(\cdot)$  is also a random variable which may not be known to the charging station. However, in Section~\ref{sec:v2g_profitability} we analyze the profit of the charging station, and the user's incentive to participate in the V2G service for a known $C_k(\cdot)$.

\subsection{Profitability of the V2G service}
Note that if $BU\geq 1$ provides a positive payoff and a  higher payoff  to the user compared to $BU=0$, then the user will opt for V2G service. The following result formalizes the condition. 
\begin{theorem}\label{thm:condn}
User $k$ opts for battery degradation $BU\geq 1$, if $u^{BU}_{k,l,t}-p^{BU}_{k,l,t}\geq \max\{0, u^{0}_{k,l,t}-p^{0}_{k,l,t}\}$ for some $BU\geq 1$. 
\end{theorem}

Now, we consider the pricing strategy stated in (\ref{eq:guranteed_profit}) i.e., $p^{BU}_{k,l,t}=v^{BU}_{l,t}-v_{-k}+\beta$. Note that with a linear  {\em battery degradation cost} discussed in the last section $u^{BU}_{k,l,t}=u^{0}_{k,l,t}-\alpha_kBU.$ Our next result characterizes the condition for user's participation in the V2G service for a linear battery degradation cost. From Theorem~\ref{thm:condn} %$L^{b}_{k,i,j}=L^{0}_{k,i,j}-\alpha b$. The following result characterizes the setting where the charging station can get strictly higher profit from the V2G service. 
\begin{theorem}
If the following is satisfied: $v^{BU}_{l,t}<v^{0}_{l,t}-\alpha_kBU$  for some $BU\geq 1$,$l$,and $t$; then the user $k$ will  have any incentive for V2G service when the pricing strategy is as given in (\ref{eq:guranteed_profit}).

The expected profit of the charging station also increases under the above condition. 
\end{theorem}
{\em Outline of the Proof}: Note that $u^{BU}_{k,l,t}-p^{BU}_{k,l,t}=u^{0}_{k,l,t}-\alpha_kBU-v^{BU}_{l,t}+v_{-k}-\beta$. Thus, the user will never select the contract with $BU=0$ if $v^{BU}_{l,t}<v^{0}_{l,t}-\alpha_kBU$ for some $BU\geq 1$. 

In the fixed profit scheme, the charging station always gets a profit of $\beta$ for if a contract is selected. When $v^{BU}_{l,t}<v^{0}_{l,t}-\alpha_kBU$ for some $BU\geq 1$, the user's probability of selecting any contract increases.  Hence, the expected profit of the charging station increases from Theorem 5. \qed

Thus, if $v^{BU}_{l,t}<v^{0}_{l,t}-\alpha_kBU$ for some $BU\geq 1$,  the user's surplus and the profit of the charging station both increase. %Note that $v^{BU}_{l,t}$ and $v^{0}_{l,t}$ can be easily computed by solving the linear programming problem. Hence, the above condition can be easily checked. 
 
If the harvested renewable energy is large enough, then the charging station can fulfill the demand using the renewable energy. Hence, the difference between $v^{BU}_{l,t}$ and $v^{0}_{l,t}$ will be not enough for the user to participate in the V2G service. On the other hand, if the harvested renewable energy is small, then the charging station may have to buy expensive conventional energy from the grid to fulfill the demand. Hence, the user will have a higher incentive to participate in the V2G service as the difference between $v^{BU}_{l,t}$ and $v^{0}_{l,t}$ may be significant. The difference is more significant when the conventional energy is more expensive (e.g. peak period). {\em Thus, the user will have a higher incentive to participate in the V2G service when the renewable energy generation is small and the cost of the conventional energy is high.}

If the storage capacity of the charging station is large, the charging station may buy energy from the grid during the off-peak period and use it during the peak period. This also reduces the difference between $v^{0}_{l,t}$ and $v^{BU}_{l,t}$. However, if the storage capacity is low, the difference between $v^{0}_{l,t}$ and $v^{BU}_{l,t}$ again increases. Hence, {\em the user's participation towards the V2G service is more likely when the storage capacity of the charging station is small and the renewable energy harvesting is low.} {\em This shows the necessity of the V2G service for better profitability as the high storage capacity is very costly to procure and the penetration of the renewable energy is still low.} %Note that if the harvested energy is not enough, then, the charging station can sell energies from the EVs during the peak hours, and charge during the other times, and can gain higher profits. Hence, it shows that when the renewable energy generation is not enough, then only users will have preference to participate in the V2G service. 

\subsection{EVs only for discharge}
So far, we assumed that EVs only come for charging where $l>0$ for menu-pricing.  However, once the V2G service proliferates, the users with their fully charged batteries  may come during the peak hours to the charging station in order to {\em only} discharge. The users can charge their batteries again in their homes  during the off-peak times. In this manner, the users can gain some profits. Though, we have not explicitly considered this scenario, our model can be easily extended to the above scenario. For example, the price $p^{BU}_{k,l,t}<0$ where $l\leq 0$  will denote that the user $k$'s EV will be discharged {\em at most} $l$ amount within deadline $t$ and additional battery utilization $BU\geq 0$. The negative price indicates that the charging station will pay to the user, as the user is delivering energy to the grid. 

Note that the {\em only} constraint needs to be changed is to replace $l$ with $|l|$ in (\ref{eq:battery_utilization}). Since $\mathcal{P}^{BU}_{l,t_{dead}}$ still remains a linear programming problem (Definition~\ref{defn:vlt}), the charging station can still find $v^{BU}_{l,t}-v_{-k}$ and can select prices according to the strategies discussed before. Hence, our results also holds in this scenario. 
\section{Numerical Results}\label{sec:simulation_results}
\subsection{Parameters and Setup}
We numerically study and compare various pricing strategies presented in this paper. We evaluate the profit of the charging station and the user's surplus  achieved in those pricing strategies. We also analyze the impact of the V2G service. 

Similar to \cite{quadratic}, the user's utility for energy $x$  is taken to be of the form 
\begin{eqnarray}\label{eq:utility}
\begin{cases}
-x^2+2rx\quad \text{if } x\leq r\nonumber\\
r^2\quad \text{otherwise}.
\end{cases}
\end{eqnarray}
Thus, the user's utility is a strictly increasing and concave function in the amount energy consumed $x$. The quadratic utility functions for EV charging have also been considered in \cite{low, low2}. Note that the user's desired level of charging is $r$. {\em We assume that $r$ is a random variable}. \cite{gov} shows that in a commercial charging station, the average amount of energy consumed per EV is $6.9$kWh with standard deviation $4.9$kWh. We assume that $r$ is a truncated Gaussian random variable with mean $6.9$kWh and standard deviation $4.9$kWh, where the truncation is to the interval $[2, 20]$.  We assume that the maximum battery capacity is $d_{max}=25$, and the minimum capacity as $d_{min}=2$. The initial battery level of a new user is assumed to be uniformly distributed in the interval $[2,25-r]$. Note that the upper bound is $25-r$, since the user's desired level of charging is $r$. %We assume $a$ is a uniform random variable in the interval $[1/20, 1/8]$. 

Following \cite{gov}, the deadline or the time spent by an electric vehicle in a commercial charging is distributed with an exponential distribution with mean $2.5$ hours. Thus, the preferred deadline ($T_{pref}$) of the user  to be an exponentially distributed random variable with mean $2.5$. Users strictly prefer a lower deadline for a given energy. We assume the utility as a convex decreasing function of the deadline. The utility of the user after the preferred deadline is assumed to be $0$. Hence, the user's utility for energy $l$ and deadline $t$ is chosen as the following
\begin{align}\label{eq:simulation}
& U_{k,l,t}=\min\{-l^2+2rl, r^2\}\times\nonumber\\
& (\exp(T_{pref}-t-t_k)-1)^{+}/(\exp(T_{pref}-t_k)-1)
\end{align}

The cost for the battery utilization is taken as $C_k(BU)=\alpha BU$. We assume that $\alpha=0.07$. Thus, the total utility for the user for contract $(l,t,BU)$ is 
\begin{align}
U^{BU}_{k,l,t}=U_{k,l,t}-\alpha BU.
\end{align}

The arrival process of electric vehicles is assumed to be a non-homogeneous Poisson arrival process since the arrival rates vary over time.  For example, during the peak-hours (8 am to 5pm) the arrival rate is higher compared to the off-peak hours. The arrival rates are chosen as $15$  ($5$, resp.) vehicles per hour during the peak period (off-peak period, resp.). We also assume that the maximum charging rate $R_{max}$  is $3.3$ kW. We also assume $R_{min}$ as $-3.3$ kW. 

The energy harvesting device is assumed to harvest energy which is a truncated Gaussian random variable with mean $2$ and standard deviation $1$ per hour.  Initial battery level is assumed to be $0$ {\em i.e.} it is fully discharged.  The prices $c_t$ for the conventional energy is assumed to be governed by Time-of-Use (ToU) time scale. Thus, the cost of buying conventional energy varies over time. We assume that the selling price, $g_t=c_t-0.001$.

\subsection{Results}
We assume that the charging station is unaware of the exact utilities. However, it knows the lower endpoints of the support sets. The charging station also knows the value of $\alpha$ and the structure of the utility function, but does not know the realized values. We consider the pricing strategy that we have introduced in Section~\ref{sec:price_uncertainty}, which is given here
\begin{eqnarray}
p^{BU}_{k,l,t}=v^{BU}_{l,t}-v_{-k}+\max_{i,j,b}\{L^{b}_{k,i,j}-v^{b}_{i,j}+v_{-k}\}^{+}+\beta.\nonumber
%p_{k,l,t}=(1+\beta)(v_{l,t}-v_{k}+\max_{i,j}\{L_{k,i,j}-v_{i,j}+v_{k}\}^{+})
\end{eqnarray}
Recall from Definition~\ref{defn:lower_endpoint} that $L_{k,l,t}$ is the lowest end-point of the utility $U_{k,l,t}.$ We study the impact of $\beta$. We also study the impact of the maximum storage capacity of the charging station $B_{max}$. Towards this end, we consider two different values of $B_{max}$; a low value, $B_{max}=5$kWh, and a high value, $B_{max}=20$kWh.

\subsubsection{Effect on Percentage of the users admitted}
Fig.~\ref{fig:betavsadmittedusers} shows that as $\beta$ increases the number of admitted users decreases.  However, the decrement is slow initially. When $\beta$ becomes  larger than a threshold, the price selected to the users becomes very large, and thus, fewer number of EVs are admitted. Note that V2G service enables the charging station to admit more users.%compared to setting where there is no V2G, more users are admitted. 

\subsubsection{Effect of $\beta$ on User's Surplus and Profit of the charging station}
\begin{figure*}
\begin{minipage}{0.23\linewidth}
\includegraphics[trim=0in 0in 0.6in 0in, width=\textwidth]{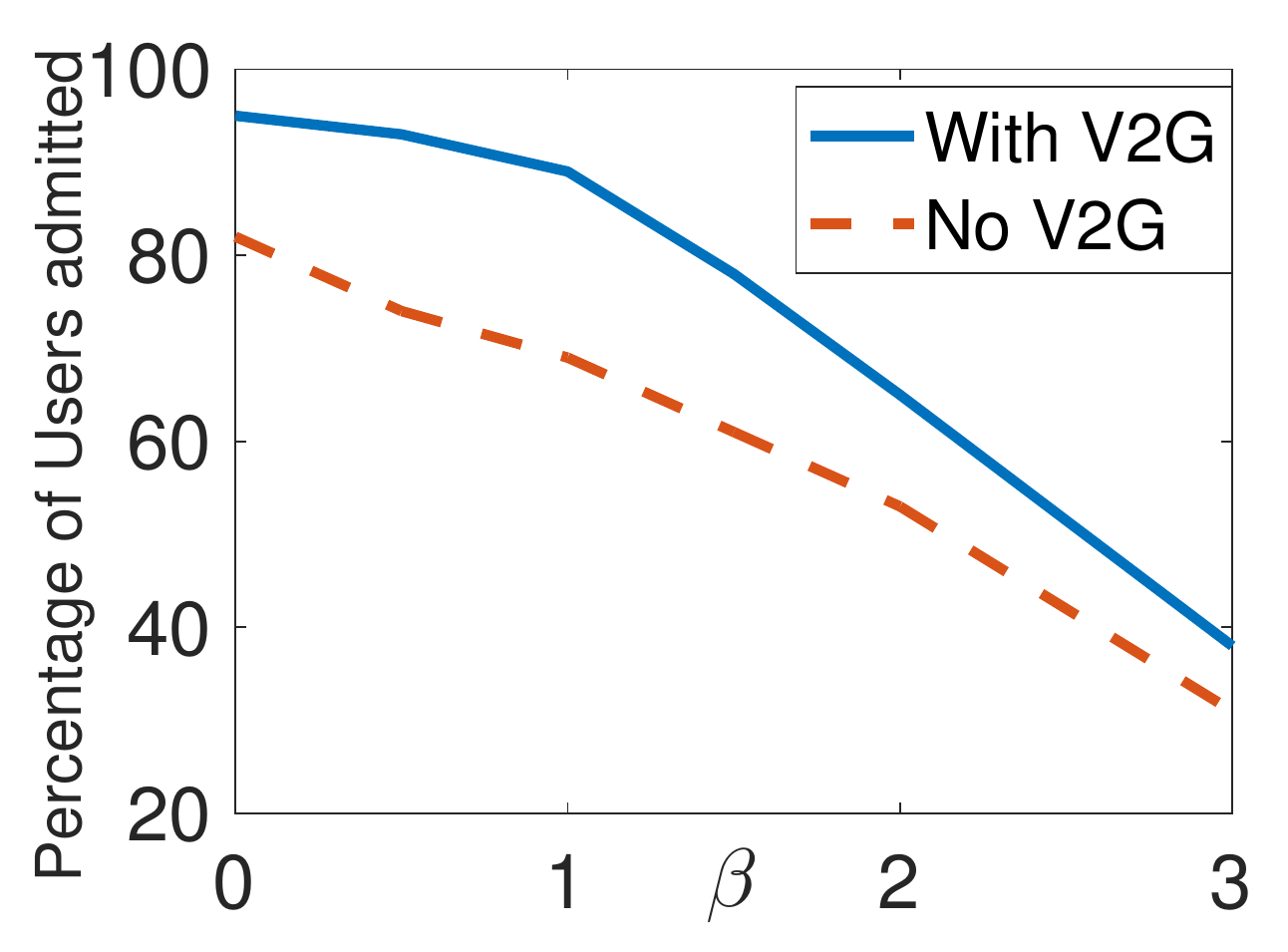}
\vspace{-0.3in}
\caption{Variation of the percentage of EVs admitted with $\beta$ with V2G service and without V2G.}
\label{fig:betavsadmittedusers}
\vspace{-0.2in}
\end{minipage}\hfill
\begin{minipage}{0.37\linewidth}
\includegraphics[trim=0in 0in 0in 0in,width=\textwidth]{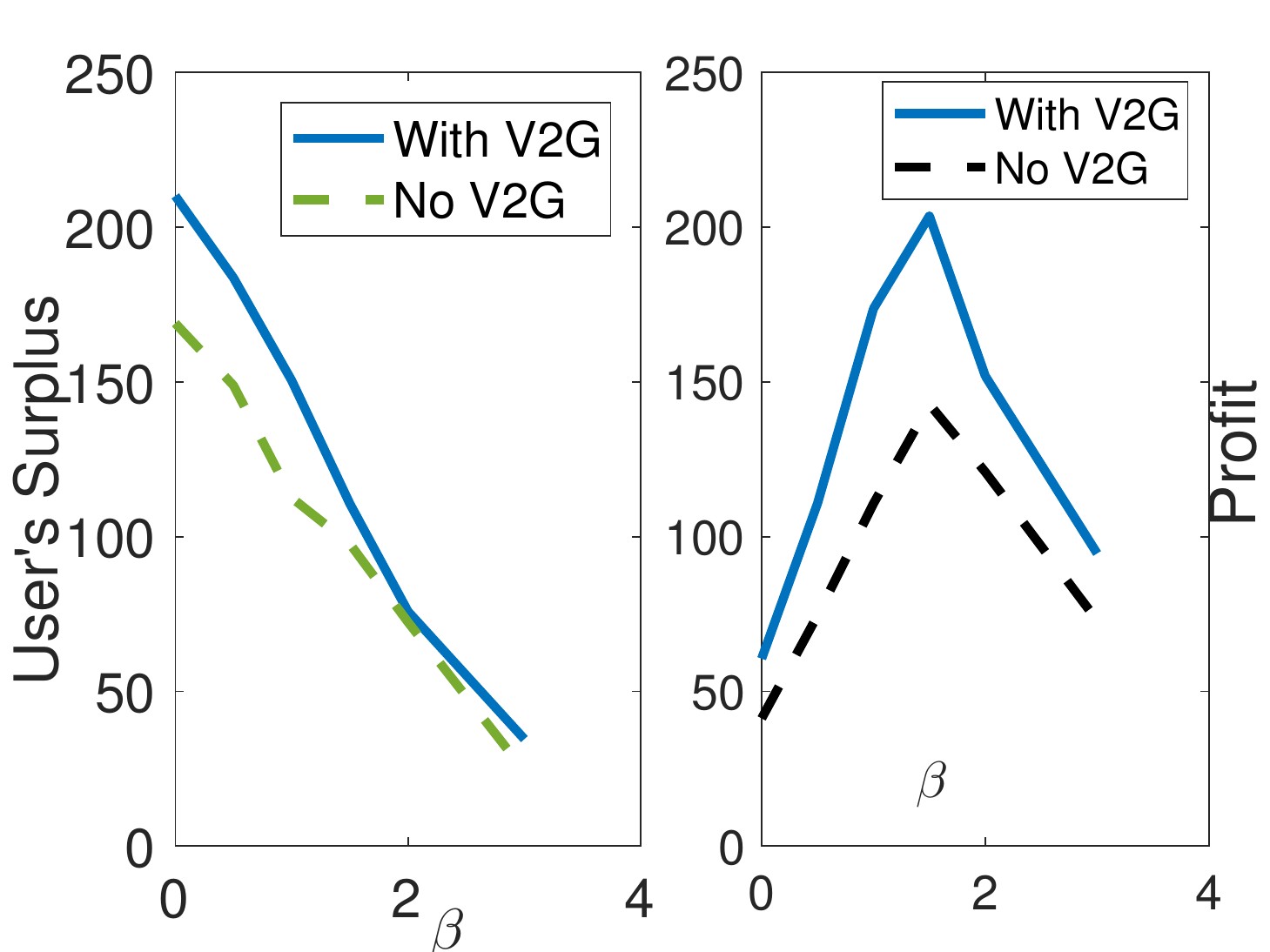}
\vspace{-0.3in}
\caption{Left-hand figure shows the variation of the total users' surpluses  with $\beta$ for $B_{max}=20$. Right-hand figure shows the profit of the station with $\beta$.}
\label{fig:consumer_surplus}
\vspace{-0.2in}
\end{minipage}
%\end{figure*}
%\begin{figure*}
\begin{minipage}{0.36\linewidth}
\includegraphics[trim=0in 0in .5in 0in,width=\textwidth]{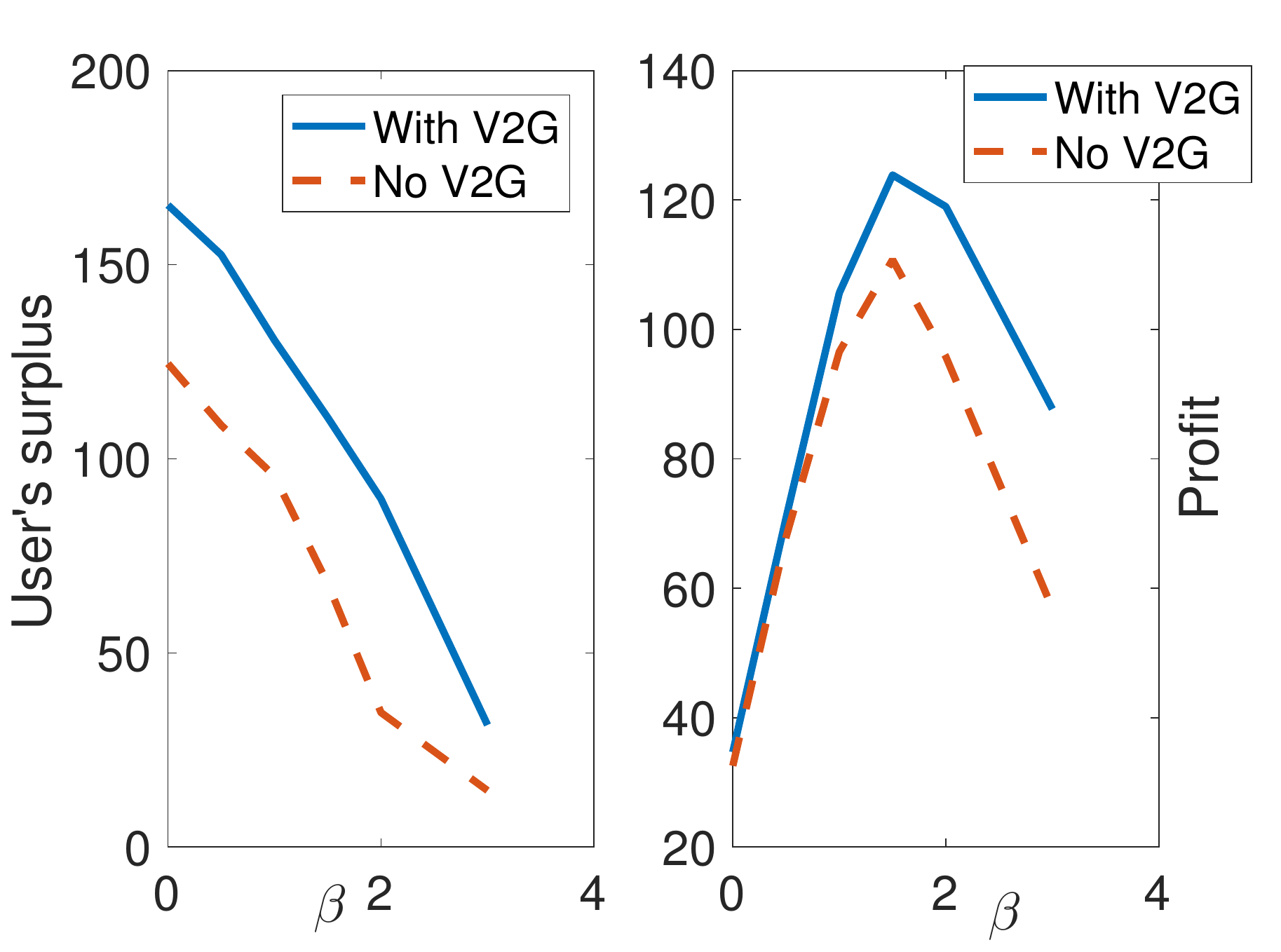}
\vspace{-0.3in}
\caption{Left-hand figure shows the variation of the total users' surpluses  with $\beta$ for $B_{max}=5$. Right-hand figure shows the profit of the station with $\beta$.}
\label{fig:soc_wellowb}
\vspace{-0.2in}
\end{minipage}
\end{figure*}

Fig.~\ref{fig:consumer_surplus} shows that as $\beta$ increases, the user's surplus decreases, as the users pay more.   The decrement is not significant for small values of $\beta$. Also, note that V2G enhances the user's surplus even though the users incur costs for the V2G service. 

As $\beta$ increases the profit increases initially. However, when $\beta>1.5$, the number of users served decreases rapidly, hence, the profit also drops.  V2G service significantly increases the profit of the charging station.

At high values of $\beta$ both users' surpluses and the profit decrease significantly. Low values of $\beta$ gives high users' surpluses, however, the profit is low.  $\beta\in [0.5,1.5]$ is the best candidate for the balance between profit and users' surpluses. 

\begin{figure*}
\begin{minipage}{0.22\linewidth}
\includegraphics[trim=0in 0in .5in 0in, width=\textwidth]{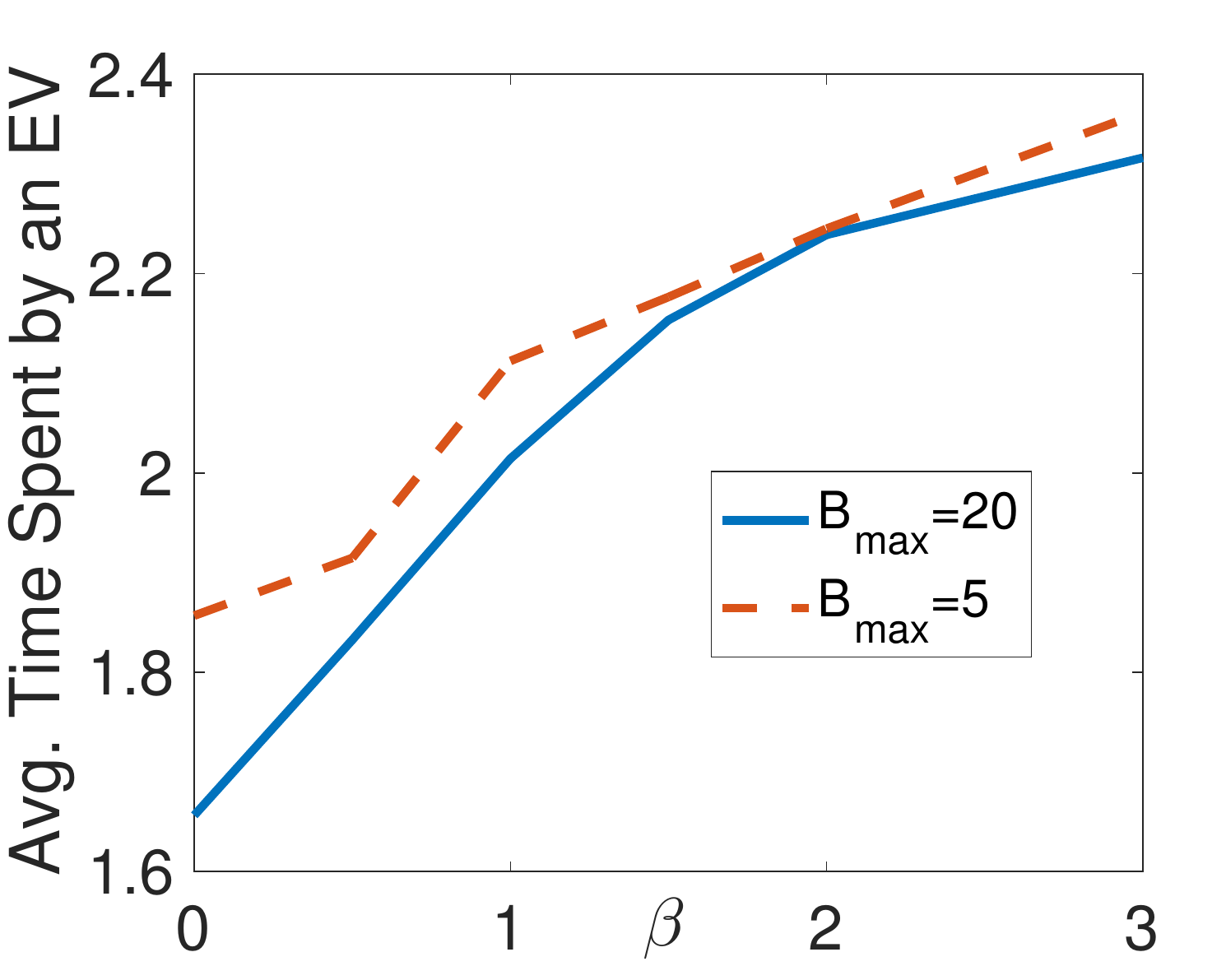}
\vspace{-0.3in}
\caption{Variation of the average time spent per EV with $\beta$.}
\label{fig:deadline}
\vspace{-0.2in}
\end{minipage}\hfill
\begin{minipage}{0.22\linewidth}
\includegraphics[trim=0in 0in .1in 0in,width=\textwidth]{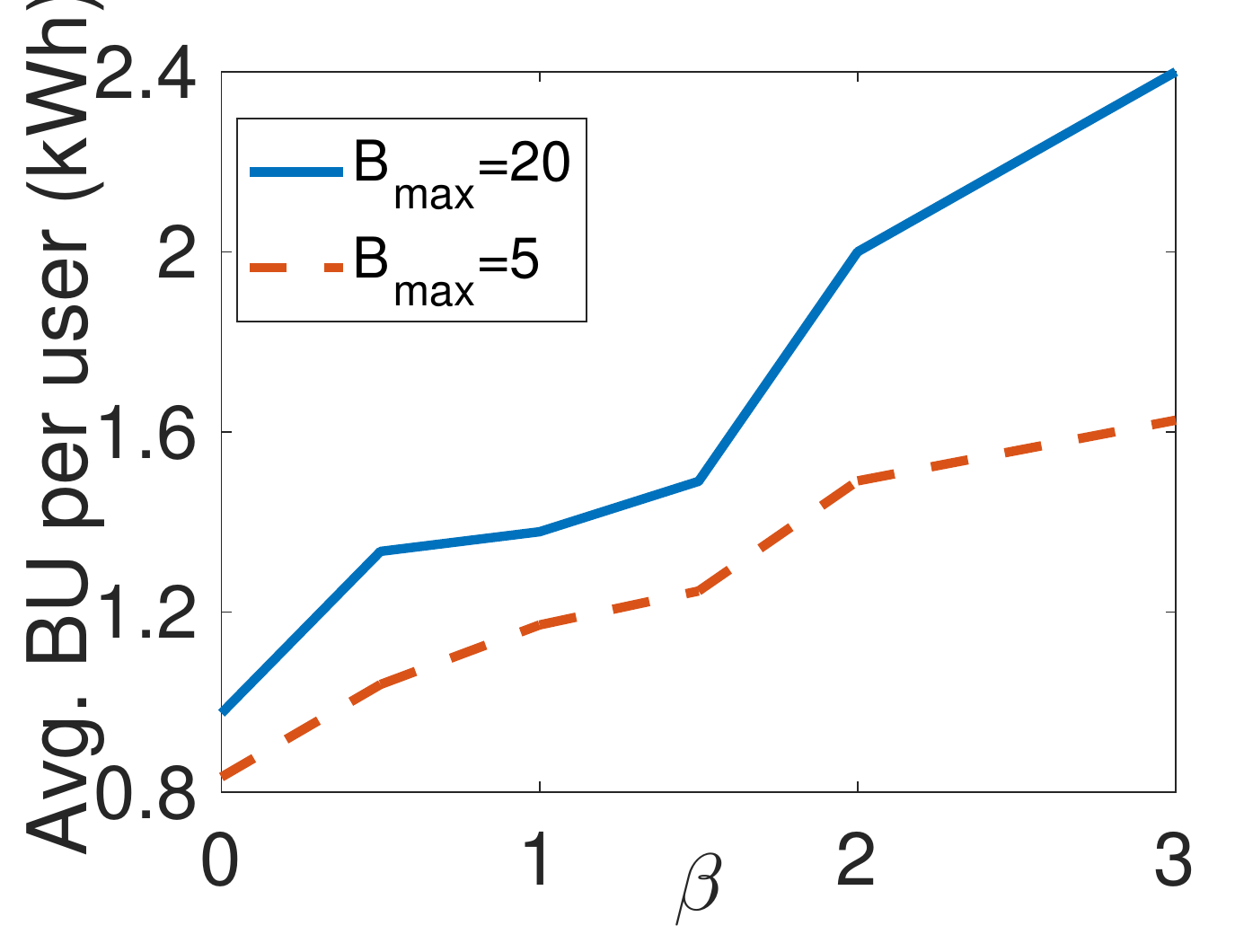}
\vspace{-0.3in}
\caption{Average amount of {\em additional} battery utilization for each EV as a function of $\beta$.}
\label{fig:bu}
\vspace{-0.2in}
\end{minipage}
\begin{minipage}{0.22\linewidth}
\includegraphics[trim=0in 0in 0.3in 0in,width=\textwidth]{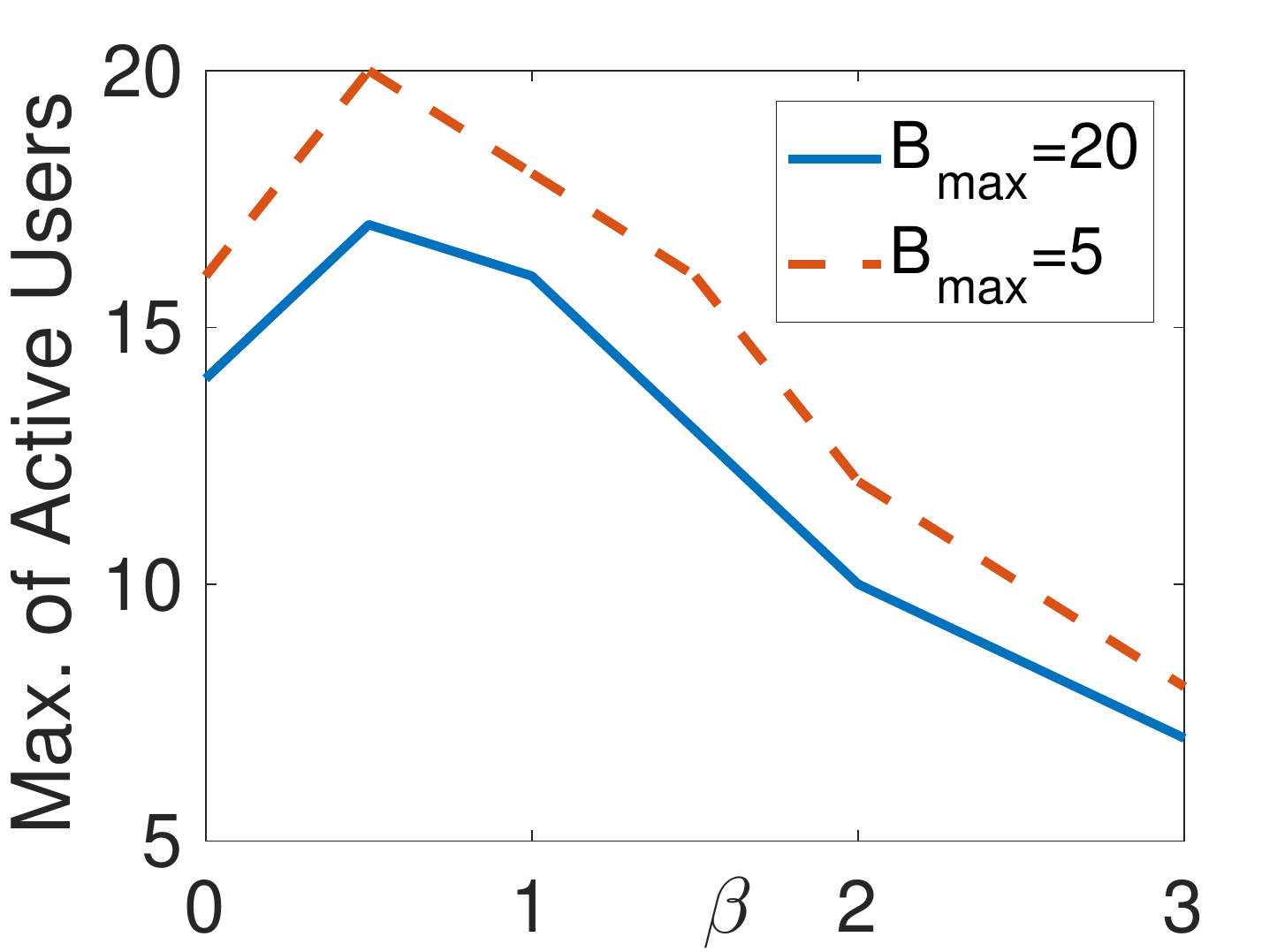}
\vspace{-0.3in}
\caption{Variation of the maximum of the number of active users  with $\beta$.}
\label{fig:active_max} 
\vspace{-0.2in}
\end{minipage}\hfill
\begin{minipage}{0.31\linewidth}
\includegraphics[trim=0in 0in .7in 0in,width=\textwidth]{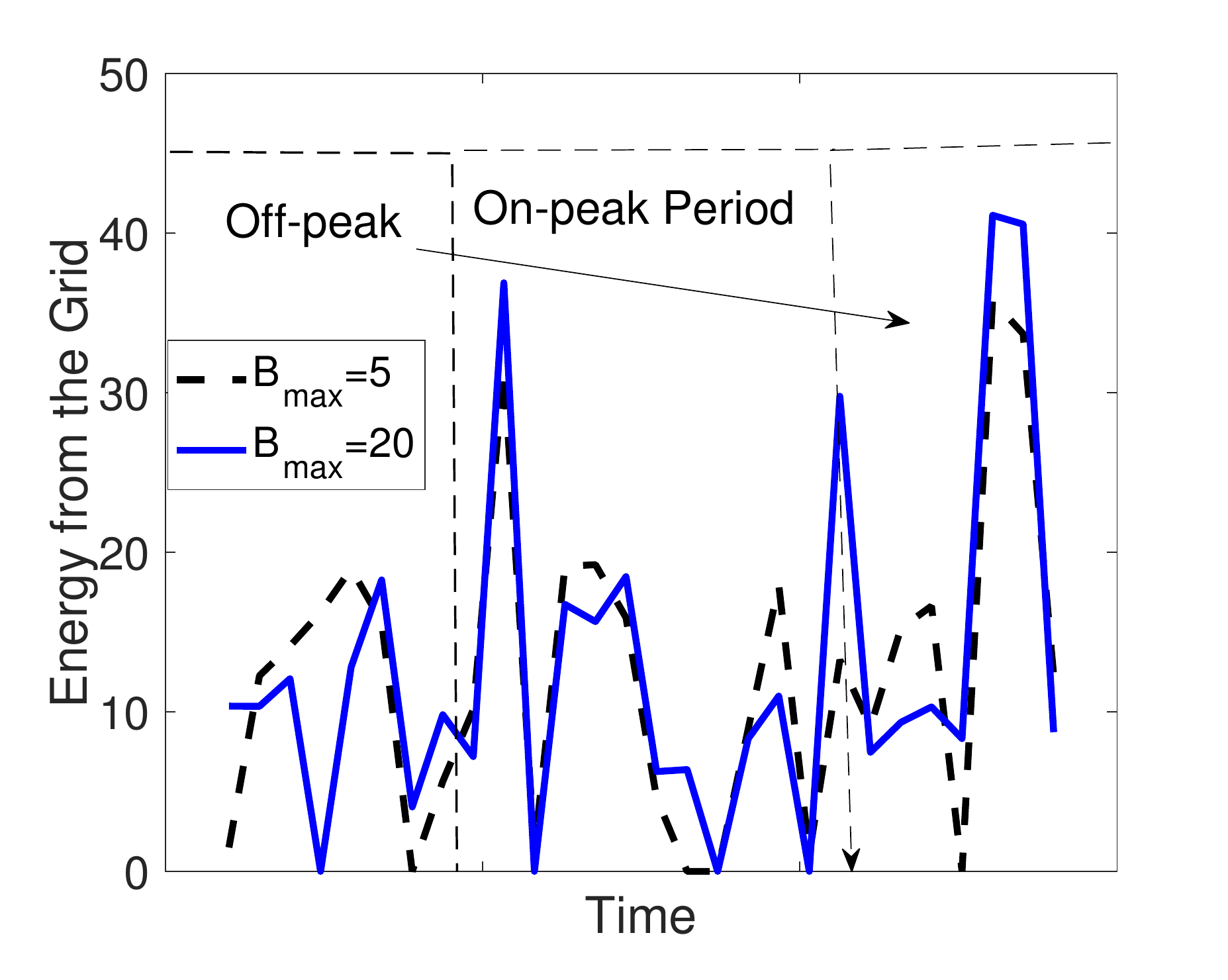}
\vspace{-0.3in}
\caption{Variation of the conventional energy bought from the grid over time for $\beta=0$.}
\label{fig:conven}
\vspace{-0.2in}
\end{minipage}
\end{figure*}
%\vspace{-0.1in}

\subsubsection{Effect of Lower $B_{max}$}
Fig.~\ref{fig:soc_wellowb} shows that if the charging station's storage capacity is low, the user's surplus and the profit decrease. This is because when the storage capacity is low, the charging station can only store lower amount of energy for the future use. Hence, the charging station has to pay more to buy conventional energy during the peak period. Hence, the charging station's profit and the user's surplus both decrease. 

%From Figs.~\ref{fig:consumer_surplus} and \ref{fig:soc_wellowb} that the V2G service fetches more profit when $B_{max}$ is low. This is because since the charging station can not store all the renewable energies generated during the off-peak hours because of the low capacity, can use the energies from the EVs during the peak hours. 

\subsubsection{Effect on the average deadline}
We subsequently study the effect of $\beta$ on the time spent by the users in the charging station (Fig.~\ref{fig:deadline}). Our analysis shows that users spend more time in the charging station with the increase in $\beta$. As $\beta$ increases, the users which have preferences for lower deadlines have to pay more, since the cost of fulfilling lower deadline contracts is  high. Hence, those users are reluctant to accept the contract. Thus, the accepted users spend more time in the charging station. 

Fig.~\ref{fig:deadline} also shows that the EVs spend lower amount of time when $B_{max}$ is high. This is because when $B_{max}$ is high, the charging station's cost to fulfill contracts with a lower deadline is lower compared to low $B_{max}$. %Though the increment of the average time spent by an EV is not exponential with $\beta$. The average time spent by an EV is  $2.5$ hours  for $\beta=1.2$ which is in accordance with the average time spent by an EV  \cite{gov}.

\begin{figure*}
\begin{minipage}{0.28\linewidth}
\includegraphics[trim=0in 0in 0.7in 0in,width=\textwidth]{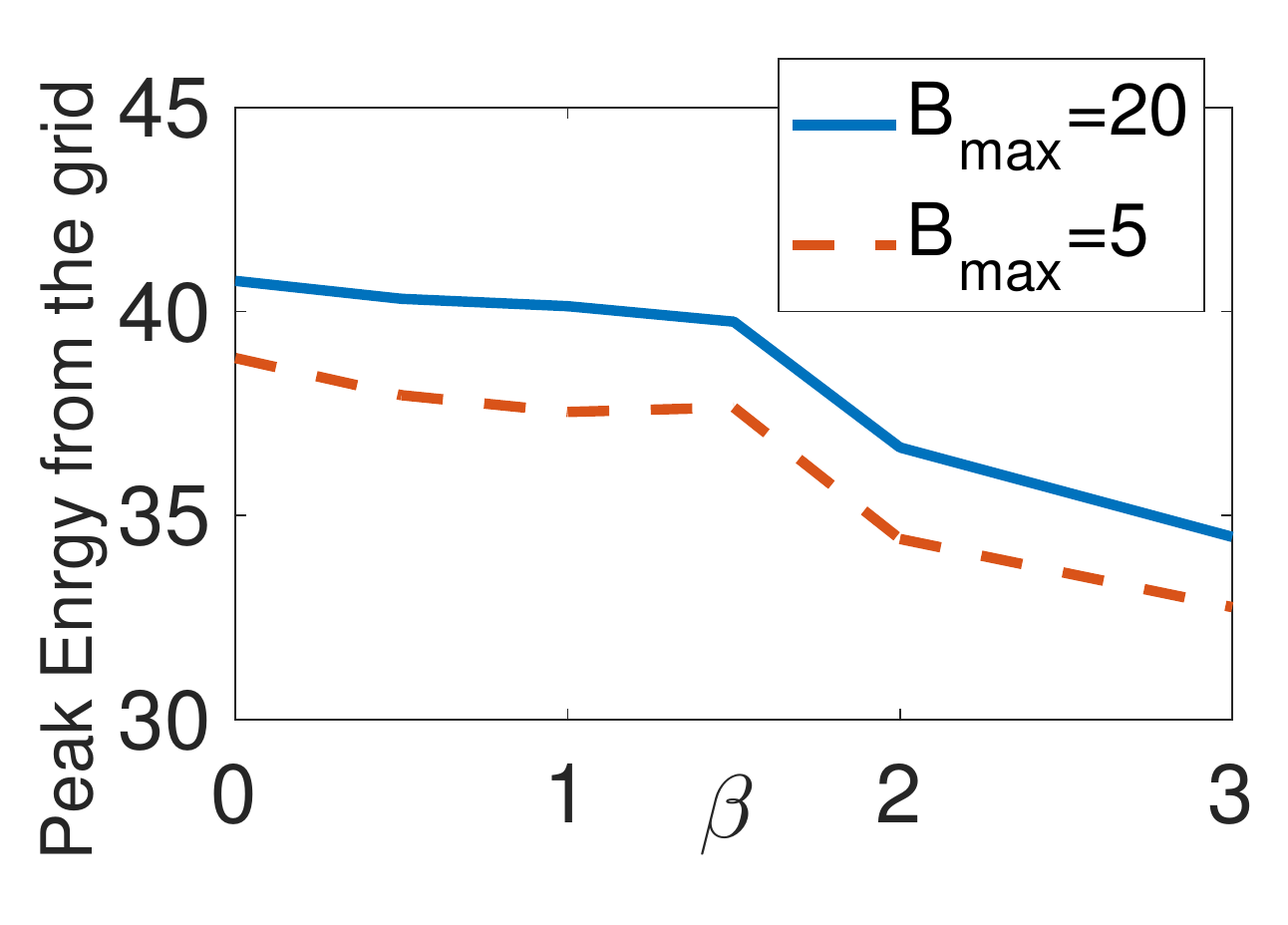}
\vspace{-0.3in}
\caption{Variation of the peak energy drawn from the grid with $\beta$.}
\label{fig:diffbuy}
\vspace{-0.2in}
\end{minipage}\hfill
\begin{minipage}{0.27\linewidth}
\includegraphics[trim=0in 0in 0.6in 0in,width=0.99\textwidth]{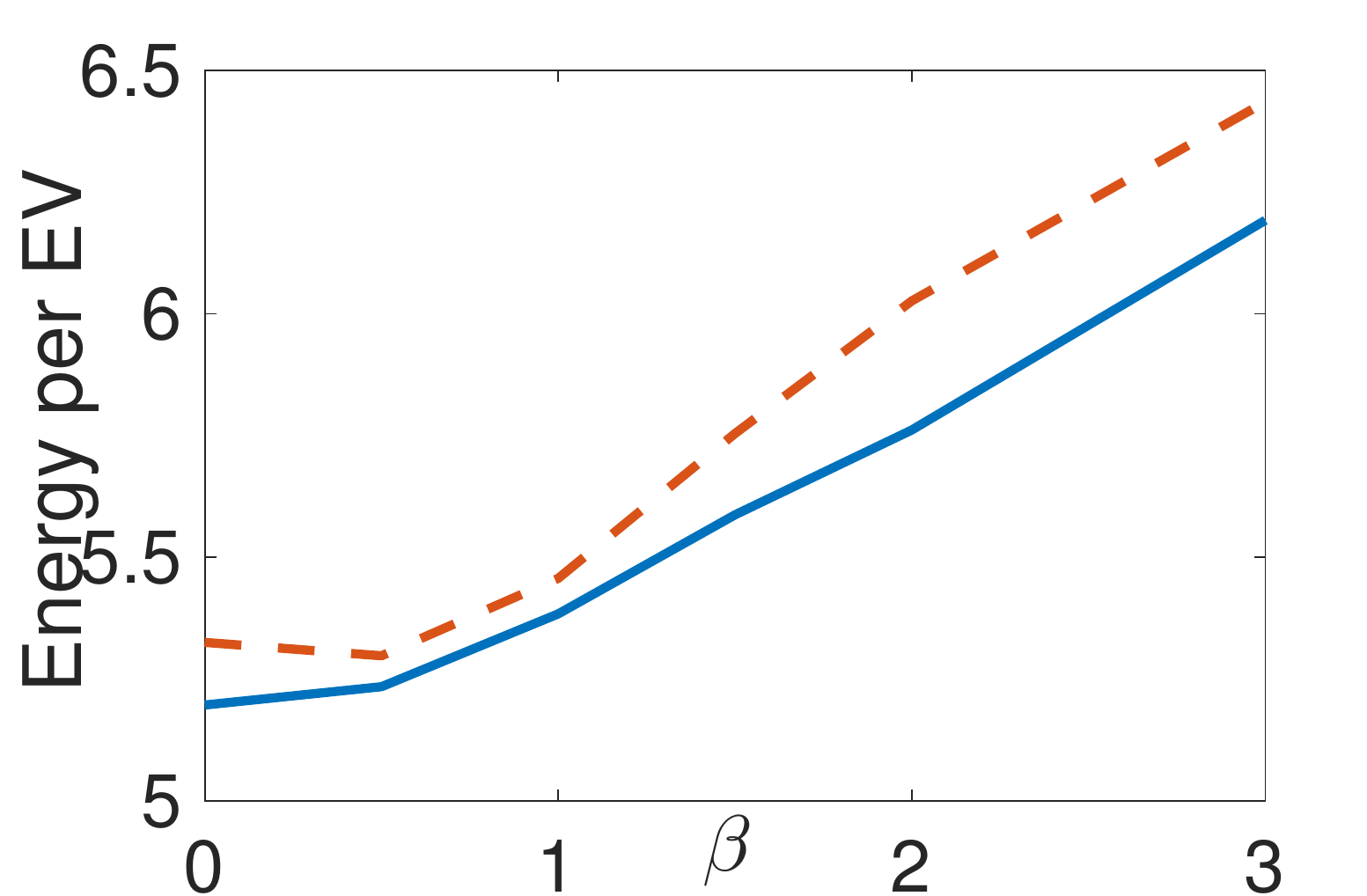}
\vspace{-0.3in}
\caption{Variation of the average energy consumed per EV with $\beta$.}
\label{fig:energy_mean}
\vspace{-0.2in}
\end{minipage}\hfill
%\end{figure*}
%\begin{figure*}
\begin{minipage}{0.33\linewidth}
\includegraphics[trim=.6in 0in 0.3in 0in,width=0.99\textwidth]{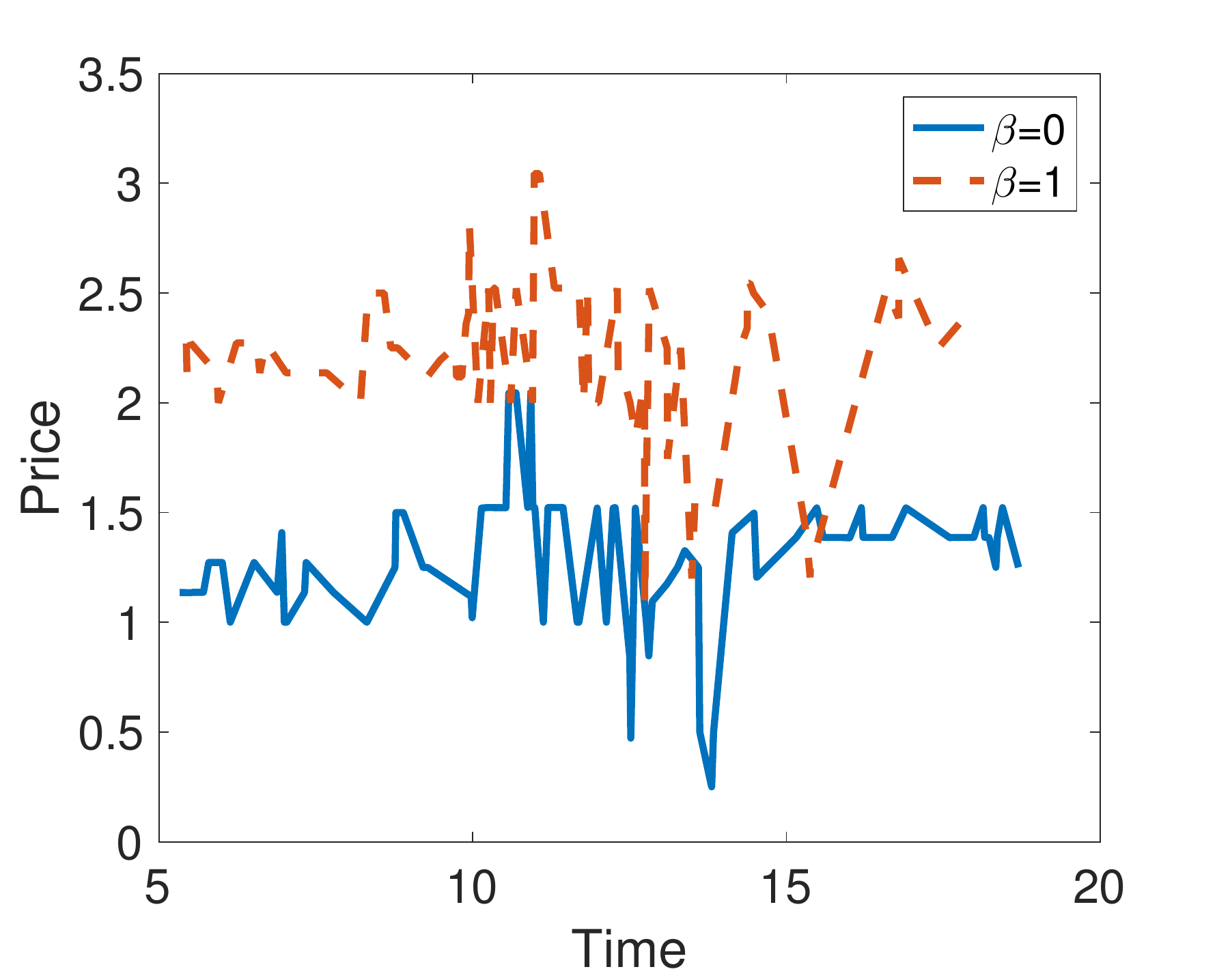}
\vspace{-0.3in}
\caption{Variation of the prices set at different times with $\beta$ for $B_{max}=20$.}
\label{fig:price}
\vspace{-0.2in}
\end{minipage}
\end{figure*}

\subsubsection{Average amount of V2G energy}
Fig.~\ref{fig:bu} shows the average amount of battery utilization of each EV as a function of $\beta$.  As $\beta$ increases the average battery utilization increases. When $\beta$ is high, the user has to pay a larger amount. However, if the user agrees for higher battery utilization, then the price may be low because the charging station may incur a lower cost. Hence, the average battery utilization increases as $\beta$ becomes high. Since the charging station can only store a lower amount of renewable energy and a lower amount of energy bought from the grid during the off-peak periods, hence, the average battery utilization is higher for lower $B_{max}$. The batteries of EVs are used more for lower storage capacity.

\subsubsection{Effect on the maximum number of active users}
Fig.~\ref{fig:active_max} shows the impact of $\beta$ on the {\em maximum} number of EVs that can simultaneously  charge (or, the number of active users) along time. This given an idea of how much charging spots we need, in order to implement the pricing strategy. 
Since the average time spent by users in the charging station increases with $\beta$ and the number of admitted users are almost the same for $\beta\leq 1$, the number of active users increases as $\beta$ increases.  Since when $B_{max}$ is low the average time spent by an EV increases, thus, the number of active users is higher when $B_{max}$ is low. %Fig.~\ref{fig:active_max} also sh%Though the maximum never reaches beyond $2$ for any value of $\beta$. However, when $\beta>1.2$, the number of active users decreases with $\beta$. 

%\subsubsection{Advantages of our proposed mechanism}
%We also show that how our pricing algorithm requires less charging spots compared to the differentiated  pricing mechanism \cite{bitar, bitar2} closest to our proposed approach (Fig.~\ref{fig:active_max}). Similar to \cite{bitar, bitar2} we assume that in a day-ahead setting the primary will set prices at every instance. The users will select the amount of energy to be consumed at by each time period. We assume that the utility will not charge beyond the preferred deadline and before the arrival time. In \cite{bitar,bitar2}  the EVs tend to spend more time as it reduces the costa and thus, the maximum of the number of the EVs present at any time is also higher (Fig.~\ref{fig:active_max}) compared to our proposed mechanism. Hence, in our proposed mechanism, the charging station controls the deadlines of the EVs and thus, lowers the number of active charging spots required and the time spent in the charging station. 

\subsubsection{Energy bought from the Grid}
Fig.~\ref{fig:conven} shows the energy bought from the grid over the day. The figure shows that the peak is significantly reduced during the on-peak period. In fact, the peak is shifted to the off-peak period even though more users come during the on peak period. This is because during the on-peak periods, the EVs are mainly used for discharging or the charging station uses more energy from its storage.  During the off-peak period, the EVs are charged, and energy is stored in the battery of the charging station. %Hence, V2G service also enables reducing the peak.the charging station has to pay more, thus, it buys more energy during the off-peak period (or, just before the beginning of the on-peak period) and stores in the battery, or the batteries of the EVs. 

\subsubsection{Impact of $B_{max}$ on the peak energy used from Grid}
Fig.~\ref{fig:diffbuy} shows that the peak energy used from the grid is higher when $B_{max}$ is high. This is because when  $B_{max}$ is high, more energy is bought and stored during the off-peak period which results in a higher peak. The reduction of the peak decreases very slowly with $\beta$ initially. However, as $\beta$ becomes very high, very few users are admitted which results in the drop of peak energy. 

\subsubsection{Effect on the average energy}
As $\beta$ increases the users with higher utilities will accept the contracts. Thus, users with higher energy demand will be more likely to accept the contract and hence, the average charging amount for each EV should increase with $\beta$. However, Fig.~\ref{fig:energy_mean} shows that for $\beta\leq 1$, the average energy consumed by each EV does not increase significantly with the increase in $\beta$. The apparent anomaly is due to the fact that  the users with higher demand but with smaller deadline preferences, may have to pay more because of the increase in the price to fulfill the contract as $\beta$ increases. Hence, such users will not accept the offers which results in a decrease of the average energy consumption with the increase in $\beta$ initially . However, as $\beta$ increases only the users with higher demand accept offers increasing the average energy consumption. Since the user spends more time for lower $B_{max}$, the average energy consumed by each EV is slightly higher for low $B_{max}$. %However the increment is only linear. In fact for $\beta=2$, the average energy consumption per EV is around $6.9$ kW-h which is in accordance with the value stated in \cite{gov}.

%\subsubsection{Effect on the Cost of the EV charging station}
%The cost of the EV charging station decreases with the increase in $\beta$ (Fig.~\ref{fig:cost}). Since the time spent by the users increases and  thus, the demand of the users can be met through renewable energies. The charging station buys lower amount of conventional energies which results into lower cost for the charging station.  When $\beta\leq 1.6$, the number of admitted users decreases sub-linearly, still the cost decreases linearly. Hence, the FERC will prefer this setting as it decreases the cost without decreasing the admitted users much.  

\subsubsection{Effect on the price selected by the charging station}
Fig.~\ref{fig:price} shows that as $\beta$ increases the price increases. Interestingly, price is not higher during the peak-period, since more EVs are involved in the V2G service during the peak period. During the off-peak period, the price from the conventional energy is low which does not provide any incentive for the V2G service.  A new price is selected every-time an EV is admitted. However, our study shows that the variation of prices within a time interval is nominal.%Hence, as $\beta$ decreases the admitted users is higher, hence the price variation is also higher as $\beta$ decreases. Also note that when $\beta$ is low, the  number of active users can be large, hence serving additional user can be significant at the peak periods as the charging station may have to buy conventional energy. Hence, the price also can be high at peak periods even when $\beta$ is low.      

%\subsubsection{Effect on the amount of battery utilization or V2G service}
%Fig.~\ref{} shows the battery utilization amount at different times. Note that during the off-peak times, the $BU$ from the users is $0$, only during the peak hours $BU$ is positive. When the battery capacity of the charging station is low, the $BU$ from the EVs is high. 

  \section{Conclusion and Future Work}
  This paper proposed a menu-based pricing approach for V2G service where the charging station selects prices to the arriving users for different charging amounts, additional amount of battery utilization, and the deadline. The user either selects one of the contracts or rejects all based on her utility function. When the charging station is aware of the utilities, there exists a profit maximizing pricing strategy which can also maximize the social welfare, however, it provides zero user's surplus. However, when the charging station is not aware of the utilities there is {\em no} profit maximizer pricing strategy which can maximize the ex-post social welfare. This paper also considered a pricing strategy which gives a fixed positive profit to the charging station irrespective of the contract chosen. The conditions where the users have incentives to participate in the V2G service are characterized, and it is seen that the V2G service is more preferable when the harvested energy is low, and the storage capacity of the charging station is low.
  
 This work can be extended in several directions. For example, the characterization of the optimal price in the non-myopic setting remain open. This paper considered a fixed price for selling the energy to the grid. The optimal price that the grid should set to enhance the V2G service or reduce the peak consumption constitutes a future research direction.

\appendix
\subsection{Proof of Theorem~\ref{thm:noposinegative}}\label{sec:proof_thmnoposinegative}
First, we show that $e_ts_t=0$ in any optimal solution of (\ref{eq:vlt}). Subsequently, we show that $r^{+}_{i,t}r^{-}_{i,t}=0$ for every $i\in\mathcal{K}_0\cup\{k\}$ and $t$. 

Suppose that in an optimal solution we have $e_t>0$ and $s_t>0$ for some $t$. We show that this can not be an optimal solution.

Consider the following
\begin{align}\label{eq:transformse}
s^{*}_t=\max\{s_t\eta_{c,cs}-e_t/\eta_{d,cs},0\},\nonumber\\
 \quad e^{*}_t=\max\{e_t/\eta_{d,cs}-s_t\eta_{c,cs},0\}.
\end{align}
The constraint in (\ref{eq:battery_capacity}) is satisfied with $e^{*}_t, s^{*}_t$ in place of $e_t$ and $s_t$ respectively. 

Now, since $\eta_{d,cs}<1$, and $\eta_{c,cs}<1$, thus, $e^{*}_t>e_t$. Thus, $e^{*}_t-s^{*}_t>e_t-s_t$. Hence,
\begin{align}\label{eq:lowervaluee}
& x^{*}_t-q^{*}_t=e^{*}_t-r_{t,\text{charge}}+r_{t,\text{discharge}}-s^{*}_t\nonumber\\
& >e_t-s_t-r_{t,\text{charge}}+r_{t,\text{discharge}}\nonumber\\
& =x_t-q_t
\end{align}
The other constraints are independent of $e_t, s_t$, hence, it will be feasible. Thus, with the modified  $x^{*}_t, q^{*}_t$, $e^{*}_t$, and $s^{*}_t$, with the other unchanged solution the value of the optimization problem is {\em lower}. Hence, this leads to a contradiction. 

Now, we show that $r^{+}_{i,t}r^{-}_{i,t}=0$ in an optimal solution of (\ref{eq:vlt}). Suppose that in an optimal solution we have $r^{+}_{i,t}>0$ and $r^{-}_{i,t}>0$ for some $i\in \mathcal{K}_0\cup \{k\}$ and $t$. Now, consider the following modification for each $i\in \mathcal{K}_0\cup \{k\}$.
\begin{align}\label{eq:transformr}
r^{+,*}_{i,t}=\max\{r^{+}_{i,t}-r^{-}_{i,t},0\}\nonumber\\
r^{-,*}_{i,t}=\max\{r^{-}_{i,t}-r^{+}_{i,t},0\}.
\end{align}
Note that because of the above modification 
\begin{align}
r^{+,*}_{i,t}-r^{-,*}_{i,t}=r^{+}_{i,t}-r^{-}_{i,t}=r_{i,t}.
\end{align}
Hence, the constraints in (\ref{eq:charged}), (\ref{eq:chargeeachev}), (\ref{eq:dischargelimit}), (\ref{eq:limit}), (\ref{eq:battery_utilization}) are satisfied for this modified solution. 

Now, consider the following 
\begin{align}\label{eq:transformrtch}
& r^{*}_{t,\text{charge}}=r^{+,*}_{k,t}/\eta_{k,c}+\sum_{i\in \mathcal{K}_0}r^{+,*}_{i,t}/\eta_{i,c},\nonumber\\
& r^{*}_{t,\text{discharge}}=\eta_{k,dc}r_{k,t}^{-,*}+\sum_{i\in \mathcal{K}_0}\eta_{i,dc}r_{i,t}^{-,*}.
\end{align}
Since $\eta_{i,c}<1$ , and $\eta_{i,dc}<1$ for all $i\in \mathcal{K}_0\cup \{k\}$ and $r^{+}_{i,t}r^{-}_{i,t}>0$ for at least index $i$, thus, $r^{*}_{t,\text{charge}}-r^{*}_{t,\text{discharge}}<r_{t,\text{charge}}-r_{t,\text{discharge}}$. Thus, if $x^{*}_t$ and $q^{*}_t$ are modified as follows
\begin{align}\label{eq:lowervaluer}
& x^{*}_t-q^{*}_t=e_t-r^{*}_{t,\text{charge}}+r^{*}_{t,\text{discharge}}-s_t\nonumber\\
& >e_t-s_t-r_{t,\text{charge}}+r_{t,\text{discharge}}\nonumber\\
& =x_t-q_t
\end{align}
then the constraints in (\ref{eq:battery_capacity}), (\ref{eq:battery_level})  are also satisfied when the rest of the solution are unchanged.   However, the cost decreases as $x^{*}_t-q^{*}_t>x_t-q_t$ with  the modified values of $r^{+,*}_{i,t}, r^{-,*}_{i,t}$, $\forall i\in\mathcal{K}_0\cup\{k\}$, $r^{*}_{t,\text{charge}}, r^{*}_{t,\text{discharge}}, x^{*}_t, q^{*}_t$. This contradicts that $r^{+}_{i,t}, r^{-}_{i,t}$ are optimal. Hence, the result follows.  \qed

The above theorem also shows the transformation of the optimal solution requires to get $e_ts_t=0$ (cf.(\ref{eq:transformse}), (\ref{eq:lowervaluee})) and $r^{+}_{i,t}r^{-}_{i,t}=0$ (cf.(\ref{eq:transformr}), (\ref{eq:transformrtch}), (\ref{eq:lowervaluer})) in case if $x_t-q_t$ are bounded. The modified solution will still be bounded.

\subsection{Proof of Theorem~\ref{thm:aclassutility}}\label{sec:proof_thmaclassutility}

Suppose the statement is false. Without loss of generality, assume that $p^{BU}_{k,l,t}=v^{BU}_{l,t}-v_{-k}+\zeta^{BU}_{l,t}$ where $\zeta^{BU}_{l,t}\neq \zeta$ for some $BU$, $l$ and $t$ achieves a strictly higher expected payoff than the pricing strategy $p^{BU}_{k,l,t}=v^{BU}_{l,t}-v_{-k}+\zeta$. 

The expected profit of the charging station for pricing strategy $p^{BU}_{k,l,t}=v^{BU}_{l,t}-v_{-k}+\zeta^{BU}_{l,t}$ is given by
\begin{align}
& \sum_{l=1}^{L}\sum_{t=t_k+1}^{T}\sum_{BU=0}^{BU_{max}}(p^{BU}_{k,l,t}-v^{BU}_{l,t}+v_{-k})\Pr(R^{BU}_{k,l,t})-v_{-k}\nonumber\\& =\sum_{l=1}^{L}\sum_{t=t_k+1}^{T}\sum_{BU=0}^{BU_{max}}\zeta^{BU}_{l,t}\Pr(R^{BU}_{l,t})-v_{-k}
\end{align}
Now, we evaluate the expression $\Pr(R^{BU}_{l,t})$. The user $k$ will select the menu $p^{BU}_{k,l,t}$ with a positive probability  if $Y^{BU}_{k,l,t}+X_{k}\geq v^{BU}_{l,t}-v_{-k}+\zeta^{BU}_{l,t}$ and for every $(i,j,b)\neq (l,t,BU)$,
\begin{align}
\zeta^{BU}_{l,t}-\zeta^{b}_{i,j}\leq Y^{BU}_{k,l,t}-v^{BU}_{l,t}-Y^{b}_{k,i,j}+v^{b}_{i,j}
\end{align}
Since $Y^{BU}_{k,l,t}, Y^{b}_{k,i,j}, v^{BU}_{l,t},v^{b}_{i,j}$ are fixed, hence,  the above inequality is either satisfied or not satisfied with probability $1$. More specifically, the user selects the menu $(l,t,BU)$ if $Y^{BU}_{k,l,t}+X_{k}\geq v^{BU}_{l,t}-v_{-k}+\zeta^{BU}_{l,t}$ and
\begin{align}
Y^{BU}_{k,l,t}-v^{BU}_{l,t}-\zeta^{BU}_{l,t}\geq \max_{i,j,b}(Y^{b}_{k,i,j}-v^{b}_{i,j}-\zeta^{b}_{i,j})
\end{align}
Without loss of generality, assume that $\zeta^{b_1}_{l_1,t_1}$ be the maximum value for which the above inequality is satisfied i.e.
\begin{align}
\zeta^{b_1}_{l_1,t_1}=\max\{\alpha^{BU}_{l,t}: Y^{BU}_{k,l,t}-v^{BU}_{l,t}-\zeta^{BU}_{l,t}\geq \max_{i,j,b}(Y^{b}_{k,i,j}-v^{b}_{i,j}-\zeta^{b}_{i,j})\}
\end{align}
The random variable $X_k$ only affects the probability whether $Y^{b_1}_{k,l_1,t_1}+X_k\geq v^{b_1}_{l_1,t_1}-v_{-k}+\zeta^{b_1}_{l_1,t_1}$ or not. Hence, the charging station's expected  profit is upper bounded by
\begin{align}
\zeta^{b_1}_{l_1,t_1}\Pr(X_k\geq v^{b_1}_{l_1,t_1}-v_{-k}+\zeta^{b_1}_{l_1,t_1}-Y^{b_1}_{k,l_1,t_1})
\end{align}
Note that by the definition of $\zeta$ (Definition~\ref{defn:alpha}),
\begin{align}\label{eq:compare}
& \zeta\max_{l,t,b}\Pr(Y^{b}_{k,l,t}+X_{k}\geq v^{b}_{l,t}-v_{-k}+\zeta)\nonumber\\
& \geq \zeta^{b_1}_{l_1,t_1}\Pr(X_k\geq v_{l_1,t_1}-v_{-k}+\zeta^{b_1}_{l_1,t_1}-Y^{b_1}_{k,l_1,t_1})
\end{align}
However by Theorem~\ref{thm:expected_payoff} the expected payoff of the charging station when it selects the price $v^{BU}_{l,t}-v_{-k}+\zeta$ is given by the expression in the left-hand side of the inequality in (\ref{eq:compare}). Hence, this leads to a contradiction. Thus, the result follows.\qed
%Consider the following two possible cases: $e_t\geq s_t$. Thus, $s^{*}_t=0$. Since $e^{*}_t<e_t$, thus, $
%
%
%
%Without loss of generality assume that
\end{document}